\documentclass[xypic,syntonly,amssymb,verbatim,12pt]{amsart}
\usepackage{graphicx}
\usepackage[all]{xy}
\usepackage{amssymb}
\usepackage{amsmath}
\usepackage[mathscr]{euscript}

\usepackage{epsfig}
\theoremstyle{plain}
\newtheorem{Thm}{Theorem}
\newtheorem{Prop}{Proposition}
\newtheorem{Cor}{Corollary}
\newtheorem{Lem}{Lemma}
\newtheorem{Def}{Definition}
 \theoremstyle{definition}
\theoremstyle{remark}
\newtheorem{Rems} {Remark}
\newtheorem{Exampl} {Example}
\numberwithin{equation}{section}
\providecommand{\norm}[1]{\lVert#1\rVert}

\begin{document}
  \title {Yang-Mills fields on $B$-branes}

 \author{ Andr\'{e}s   Vi\~{n}a} 
\address{Departamento de F\'{i}sica. Universidad de Oviedo.    Garc\'{\i}a Lorca 18.
     33007 Oviedo. Spain. }
\email{vina@uniovi.es}

 \maketitle
\begin{abstract}
Considering the $B$-branes over a complex manifold $Y$ as objects of the bounded derived category $D^b(Y)$, we define holomorphic  gauge fields on $B$-branes    and the Yang-Mills functional for these fields.
These definitions are a generalization to $B$-branes  of concepts that are well known in the context of vector bundles. 
Given ${\mathscr F}^{\bullet}\in D^b(Y)$, we show that  the Atiyah class $a({\mathscr F}^{\bullet})
\in{\rm Ext}^1({\mathscr F}^{\bullet},\,\Omega^1({\mathscr F}^{\bullet}))$ is the obstruction to the existence of gauge
 fields on ${\mathscr F}^{\bullet}$. We determine the $B$-branes over $\mathbb{ CP}^n$ that admit  holomorphic gauge fields.
We  prove that the set of Yang-Mills fields on the $B$-brane  
 ${\mathscr F}^{\bullet} $, if it is nonempty, is in bijective correspondence with the points of an algebraic subset of ${\mathbb C}^m$ defined by $m\cdot s$ polynomial equations of degree $\leq 3$, where $m={\rm dim}\,{\rm Hom}({\mathscr F}^{\bullet},\,\Omega^1({\mathscr F}^{\bullet}))$ and $s$ is the number of non-zero cohomology sheaves ${\mathscr H}^i({\mathscr F}^{\bullet})$.
 We show sufficient  conditions under them any Yang-Mills field on a reflexive sheaf of rank $1$ 
is flat.

\end{abstract}


\section {Introduction} \label{S:intro}

 In this article, we extend  the well-known concepts of gauge field and Yang-Mills field on vector bundles to $B$-branes.
From the mathematical point of view, a $B$-brane over  a compact connected complex $n$-manifold $Y$ is an object of $D^b(Y)$, the bounded derived category of  coherent analytic sheaves over $Y$ \cite[Sect. 5.4]{Aspin} \cite[Sect. 5.3]{Aspin-et}.

A holomorphic vector bundle $V$ over $Y$ is a particular case of $B$-brane over $Y$.
  A holomorphic gauge field  on $V,$ in mathematical terms a holomorphic connection on $V$ \cite{Atiyah}, allows us
 to define a derivative of the holomorphic sections of $V$ along any ``direction'' in $Y,$ giving rise to {\em holomorphic} sections. That is, the holomorphic gauge field gives identifications between the ``infinitesimally close'' fibers of $V.$  Conversely, such consistent identifications determine a gauge field.
 
 Not every holomorphic vector bundle  supports holomorphic connections,  
 unlike what happens in the smooth category;  for example, the vanihing of the Chern class is a necessary and sufficient condition for the existence holomorphic connections on a line bundle (see Appendix).  Moreover,  if the set of these gauge fields is non-empty, it is an affine space associated to a finite dimensional vector space. These properties will hold in the extension of the holomorphic gauge fields to more general $B$-branes.


 \smallskip
{\it Gauge fields on coherent sheaves.}	
	Given  a holomorphic vector bundle $V$ over the manifold $Y,$ a holomorphic gauge field 
	on $V$ can be regarded as a right inverse of the projection $\pi:J^1(V)\to V,$ where $J^1(V)$ is the $1$-jet bundle of $V$. For details see \cite[Sect. 3.1]{Vina21}.

	That approach  admits a natural translation to the context of the coherent sheaves. If ${\mathscr F}$ is a coherent sheaf over the compact analytic manifold $Y,$ we define a holomorphic gauge field on ${\mathscr F}$ as a right inverse of the natural morphism ${\mathscr J}^1({\mathscr F})\to{\mathscr F},$ where ${\mathscr J}^1({\mathscr F})$  is the corresponding  $1$-jet sheaf. Denoting by $\Omega^p$ the sheaf of the holomorphic $p$-forms on $Y,$ that inverse determines a morphism of {\it abelian} sheaves  
	$$\nabla:{\mathscr F}\to \Omega^1({\mathscr F}):=\Omega^1\otimes_{{\mathscr O}}{\mathscr F},$$ which satisfies the Leibniz's rule.  Conversely, such a morphism defines a holomorphic gauge field in the above sense.
	
	The obstruction to the existence of a holomorphic connection on the sheaf ${\mathscr F}$ is 
	 an element of the group ${\rm Ext}^1({\mathscr F},\,\Omega^1({\mathscr F}))$. Furthermore, when the set of holomorphic gauge 
	 fields on ${\mathscr F}$ is nonempty, it is an  affine space   associated  to the {\it finite} dimensional vector space  
	 ${\rm Hom}({\mathscr F},\,\Omega^1({\mathscr F}))$.
 
 The curvature of the holomorphic gauge field $\nabla$  
 is defined in the usual way and turns out to be an element ${\mathcal K}_{\nabla}\in {\rm Hom}({\mathscr F},\,\Omega^2({\mathscr F)})$. When ${\mathcal K}_{\nabla}=0$, we say that $\nabla$ is {\it flat}.

\smallskip	

	If $\nabla$ is a  flat  holomorphic connection on ${\mathscr F}$, then one has the complexes
	$\big({\mathscr A}^{\bullet, 0}({\mathscr F}):={\mathscr A}^{\bullet, 0}\otimes_{\mathscr O}{\mathscr F},\, \nabla\big)$,
	where ${\mathscr A}^{p, q}$ is the sheaf of $C^{\infty}$ differential forms of type $(p,q)$ on $Y,$  the complex $(\Omega^{\bullet}({\mathscr F}),\,\nabla)$ and
	 $({\mathscr A}^{\bullet}({\mathscr F}),\, \nabla+\bar\partial).$  
	
	If on every fibre ${\mathscr F}_{(x)}$ is defined an Hermitian metric $\langle\,,\,\rangle_x$ such that 
	$\{\langle\,,\,\rangle_x\}_x$ is a ``smooth'' family  (see Definition \ref{D:HermitianSheaf}), we say that ${\mathscr F}$ is a Hermitian sheaf. When ${\mathscr F}$ is a locally free sheaf, this definition coincides with the usual on holomorphic vector bundles \cite[Chap. III]{Wells}.
	
 Using the index formula, we prove the following theorem, which   asserts that the Euler characteristic of the above complexes, when 
 ${\mathscr F}$ is a locally free Hermitian sheaf, is determined by certain characteristic classes of $Y$ and the rank of ${\mathscr F}$.  
 \begin{Thm}\label{IndexTheorem}
 Let ${\mathscr F}$ be a Hermitian locally free sheaf   of rank $r$ over the K\"ahler $n$-manifold $Y$. If $\nabla$ is a holomorphic flat gauge field on ${\mathscr F}$, then 
$$\chi(\Omega^{\bullet}({\mathscr F}))=r\,{\rm e}(Y)=\chi({\mathscr A}^{\bullet}({\mathscr F})),\;\;\;
\chi({\mathscr A}^{\bullet, 0}({\mathscr F}))=(-1)^nr\,{\rm td}_{\mathbb C}(\bar Y),$$
Where $\chi({\mathscr C}^{\bullet})=\sum_i(-1)^i{\rm dim}\,\Gamma(Y,{\mathscr C}^i)$, ${\rm e}(Y)$ is the Euler class of $Y,$
and 
${\rm td}_{\mathbb C}(\bar Y)$ is the complexified Todd class of manifold conjugated of $Y$ (i.e., if the complex estructure of $Y$ is defined by the endomorphism $J$, the one of $\bar Y$ is defined by $-J$).
 \end{Thm}


	\smallskip
	{\it Yang-Mills fields on shaves.}
If $Y$ is a K\"ahler manifold, on  the set of holomorphic connections over the Hermitian sheaf ${\mathscr F}$, one defines the Yang-Mills functional $\mathcal{YM}$. The value  of  $\mathcal{YM}$ at a connection $\nabla$
	is the squared norm 
	$\norm{{\mathcal K}_{\nabla}}^2$ of its curvature.
	 The stationary points of functional $\mathcal{YM}$ are the Yang-Mills gauge fields on ${\mathscr F}$. The set of such critical points
	 will be denoted by ${\sf YM}({\mathscr F})$.
		We will prove the following theorems.
	\begin{Thm}\label{P:numeroYM}
	If ${\mathscr F}$  admits a holomorphic gauge field and 
	$$m={\rm dim}\, {\rm Hom}({\mathscr F},\,\Omega^1({\mathscr F}) ).$$ 
	Then the set ${\sf YM}({\mathscr F}),$ of holomorphic Yang-Mills fields on ${\mathscr F},$ is in bijective correspondence with the points of an algebraic set in ${\mathbb C}^m$ defined by $m$   algebraic equations of degree $\leq 3$. In particular, if 
	 $m=2$ and the cardinal of ${\sf YM}({\mathscr F})$ is finite, then
	$\# {\sf YM}({\mathscr F})  \leq 9$. 
	\end{Thm}
	
	\begin{Thm}\label{impatial=0}
	Let ${\mathscr F}$ be a coherent reflexive  Hermitian  sheaf with rank $1$ over  a Hodge manifold,
	such that $c_1({\mathscr F} )=0$. Then a holomorphic gauge field on   ${\mathscr F}$ is Yang-Mills, iff it is flat.  
	\end{Thm}
	

\smallskip

{\it Yang-Mills fields on a $B$-brane.}	
	As we said, a $B$-brane on $Y$ is 
	a complex   $({\mathscr F}^{\bullet},\,\delta^{\bullet})$ of analytic coherent sheaves on $Y$. According to the preceding paragraphs, it is reasonable to define a gauge field on this brane as an element of
 ${\rm Hom}_{D^b(Y)}\big({\mathscr F}^{\bullet},\,  {\mathscr J}^1({\mathscr F}^{\bullet})\big)$ which lifts the identity on  ${\mathscr F}^{\bullet}$.
 
 Using that the derived category $D^b({\mathbb P}^n)$ is generated by  the family
 \begin{equation}\label{set}
{\sf E}=\{{\mathscr  O}_{{\mathbb P}^n}(-n),\, {\mathscr  O}_{{\mathbb P}^n}(- n +1),\dots, {\mathscr O}_{{\mathbb P}^n}(-1),\,{\mathscr O}_{{\mathbb P}^n}\},
\end{equation} 
 we  prove the following theorems. 
 \begin{Thm}\label{Th:card} 
 The cardinal of the set of holomorphic gauge fields on any $B$-brane over ${\mathbb P}^n$ is $\leq 1.$
 \end{Thm}

Particular $B$-branes over ${\mathbb P}^n$ are the  complexes consisting of direct sum of copies of ${\mathscr O}_{{\mathbb P}^n}$
\begin{equation}\label{bigoplusimp1}
 \dots \to \bigoplus_{i\in S_p}{\mathscr O}_{{\mathbb P}^n}\overset{d^p}{\longrightarrow} \bigoplus_{i\in S_{p+1}}{\mathscr O}_{{\mathbb P}^n}\to\dots
 \end{equation}
where the $S_p$ are finite sets and the coboundary operators are constant matrices. We will prove the following theorem.
\begin{Thm}\label{Th:2}
 A $B$-brane on ${\mathbb P}^n$ admits a holomorphic gauge field iff it is isomorphic to a brane 
  of the form (\ref{bigoplusimp1}).
 \end{Thm}

  When $Y$ is a Hodge manifold, it is a projective smooth variety and by the GAGA correspondence the analytic coherent sheaves on $Y$ can be considered as algebraic ones. 
  Moreover, the category $D^b(Y)$ is equivalent to the homotopy category $K^b(Y)_{\rm coh}$ of complexes of injective sheaves with  coherent and bounded    cohomology   (see Subsection \ref{homotopycategory}).
   Considering ${\mathscr F}^{\bullet}$ as an object of $K^b(Y)_{\rm coh}$, we set $\widehat{\Omega^1({\mathscr F}^{\bullet})}$ for an object of $K^b(Y)_{\rm coh}$ quasi-isomorphic to $\Omega^1({\mathscr F}^{\bullet})$.
   In this way, a gauge field $\psi$ on ${\mathscr F}^{\bullet}$ is a 
    homotopy class of morphisms between the complexes  ${\mathscr F}^{\bullet}$ and $\widehat{\Omega^1({\mathscr F}^{\bullet})},$ considered as complex of abelian sheaves.
	Thus, that homotopy class determines for each $j$ a unique morphism
	$\vartheta^j:{\mathscr H}^j({\mathscr F}^{\bullet})\to {\mathscr H}^j(\Omega^1({\mathscr F}^{\bullet})),$
	 between   the corresponding cohomology sheaves.
  In summary, the gauge field on ${\mathscr F}^{\bullet}$ defines a family of connections $\vartheta^j$ on the cohomology sheaves
  ${\mathscr H}^j({\mathscr F}^{\bullet})$.

 When the cohomology sheaves ${\mathscr H}^j$ are Hermitian, in which case we say that ${\mathscr F}^{\bullet}$ is Hermitian, we define the value of the Yang-Mills functional on the above  gauge $\psi$ as
 $\sum_i(-1)^i\norm{{\mathcal K}_{\vartheta^i}}^2.$ Thus, the Yang-Mills functional is a kind of Euler characteristic  of the gauge field. 
  This definition of the Yang-Mills functional on branes  generalizes the one given for coherent sheaves, obviously.

The gauge field $\psi$ is a Yang-Mills field if it is a stationary point of Yang-Mills functional. 
We will also prove the following result.
 
\begin{Thm}\label{C:definitivo}
Let $({\mathscr F}^{\bullet},\,\delta^{\bullet})$ be a Hermitian  $B$-brane on the Hodge manifold $Y,$ such that sheaves ${\mathscr H}^i({\mathscr F}^{\bullet})$ are reflexive. A gauge field 
$\psi$ on the brane is   a Yang-Mills field iff  $\vartheta^i$ is a Yang-Mills field
	on ${\mathscr H}^i({\mathscr F}^{\bullet})$ for all $i$.
\end{Thm}
 In Theorem \ref{P:numero YM}, we generalize the result given in Theorem \ref{P:numeroYM}, about the cardinal of the set of Yang-Mills fields on a sheaf, to a general brane.
 \smallskip
 

The article is organized in two sections. In Section \ref{S:YMandF}, we define the gauge fields on a coherent sheaf.
 In the first subsections of this section, we revise the definition of the $1$-jet sheaf of a coherent sheaf ${\mathscr F}$
 because,
 although this is well known in algebraic geometry, it is not so well known in the community of mathematical physicists.
 We prove also Theorem \ref{IndexTheorem}, above mentioned.
  The Yang-Mills functional is introduced in Subsection \ref{ss:holomorphicYang}, where Theorems 
 \ref{P:numeroYM} and \ref{impatial=0} are proved. We also describe some properties of the Yang-Mills fields on reflexive sheaves.
  

 In Section \ref{S:Fields_brane} are considered the gauge fields on a general $B$-brane. In Subsection \ref{Ss:Fields_brane},
 we give the definition of gauge field on a $B$-brane and prove  Theorems \ref{Th:card} and \ref{Th:2} about the existence of gauge fields on branes over ${\mathbb P}^n.$   We define in Subsection \ref{SsYng-Millsonabrane} the Yang-Mills   functional, showing the reasons on which this definition is based and prove  
 Theorems  \ref{C:definitivo} and \ref{P:numero YM}.
   
    For the sake of completeness we give in the Appendix a simple proof,
  in the context of the $\check{{\rm C}}$ech cohomology, of 
  the following well-known fact in algebraic geometry and which has been mentioned above:
  A necessary and sufficient condition for the existence of
    holomorphic connections on a line bundle $L$ is the vanishing of $c_1(L)$. This result can certainly be generalised to  bundles of arbitrary rank \cite{Atiyah}, but perhaps such a simple proof might be of interest to some mathematical physicists.


\section{Gauge fields on coherent sheaves}\label{S:YMandF}


\subsection{Gauge fields on a sheaf.} \label{ss:FirstJet} 
As explained in the Introduction,
 the existence of a gauge field on a coherent sheaf ${\mathscr F}$ over a  compact connected complex  manifold $Y$ should define an isomorphism between the stalks  of ${\mathscr F}$ at any ``infinitesimally close'' points of $Y$. 
 
 The idea of being infinitesimally close can formulated by means of the first infinitesimal neighborhood $Y^{(1)}$ of the diagonal of $Y$ \cite[page 698]{G-H}. 
 If $R$ is a ${\mathbb C}$-algebra, ${\rm Hom}({\rm Spec}\, R,\,Y)$ is the set of points of $Y$ with values  in $R$. Two points $x_1,\,x_2$ are infinitesimally close if the morphism $(x_1,\,x_2):{\rm Spec}\, R\to Y\times Y$ factorizes through $Y^{(1)}$ \cite[page 6]{Deligne}. In this case, there is a morphism
 $h:{\rm Spec}\,R \to Y^{(1)}$ such that  $\pi_i h=x_i$,  where $Y\overset{\pi_1}{\leftarrow}Y^{(1)}\overset{\pi_2}{\to} Y$  are the projections. Given a sheaf ${\mathscr F}$ on $Y,$ each element $\alpha$ of the set 
\begin{equation}\label{Deligne_Gauge}
 {\rm Hom}_{Y^{(1)}}(\pi_1^*{\mathscr F},\, \pi_2^*{\mathscr F}),
 \end{equation}
gives rise to the morphism $h^*(\alpha):x_1^*{\mathscr F}\to x_2^*{\mathscr F}.$
 Therefore, following Deligne, one can consider a  holomorphic gauge field on ${\mathscr F}$ as an element of  (\ref{Deligne_Gauge})
 which is the identity on $Y.$
 
 On the other hand, each element of (\ref{Deligne_Gauge})
determines,  via the adjunction isomorphism, a morphism ${\mathscr F} \to\pi_{1*} \pi^*_2{\mathscr F},$ and conversely.

As  $\pi_{1*} \pi^*_2{\mathscr F}$ is
 the first jet sheaf ${\mathscr J}^1({\mathscr F} )$ of the coherent sheaf ${\mathscr F},$
in the following paragraphs we review the definition of the  jet sheaf and also that of first neighborhood of the diagonal.


 \subsubsection{Neighborhood of the diagonal.}\label{Ss:Neighborhood}

 We denote by $i:\Delta\hookrightarrow Y\times Y$ the embedding of the diagonal. As a closed subvariety, $\Delta$ is defined by an ideal $\mathscr I$ of
 $\Hat{\mathscr O}:={\mathscr O}_{Y\times Y}.$ We will consider the following ringed spaces
$$(Y,\,{\mathscr O}),\;\;Y^{(1)}=\Big(\Delta,\,{\mathscr O}_{Y^{(1)}}:=\big(\Hat{\mathscr O}/{\mathscr I}^2\big)|_{\Delta}\Big),\;\; (Y\times Y,\,\Hat{\mathscr O}).$$
 $Y^{(1)}$ is the first infinitesimal neighborhood of $\Delta$. 

\smallskip
\noindent
{\it The ${\mathscr O}\mbox{-}{\mathscr O}$-bimodule structure in ${\mathscr O}_{Y^ {(1)}}$}.
 We set $p_1,p_2: Y\times Y\rightrightarrows Y$ for the corresponding projection morphisms. The natural morphisms between the above topological spaces are shown in the following commutative diagram
 \begin{equation}
  \xymatrix{
 Y\ar[r]_{k} \ar@/^1pc/[rr]^j &\Delta \ar[r]_i &
 Y\times Y\ar@/^/[r]^{p_1}\ar@/_/[r]_{p_2}    & Y
 }
 \end{equation}

Given $f\in {\mathscr O}$ and $l\in\Hat{\mathscr O}$, the product
\begin{equation}\label{right-product}
(l+{\mathscr I}^2)\cdot f=l\,(f\circ p_2)+ {\mathscr I}^2
\end{equation}
defines a right ${\mathscr O}$-module structure on $\Hat{\mathscr O}/{\mathscr I}^2$.
More explicitly,  $(l(x,\,y)+{\mathscr I}^2)\cdot f=l(x,\,y)f(y)+ {\mathscr I}^2.$
Analogously, $f\cdot(l+{\mathscr I}^2)=(f\circ p_1)l+ {\mathscr I}^2$ gives to $\Hat{\mathscr O}/{\mathscr I}^2$
a left ${\mathscr O}$-module structure. However,  the restrictions to ${\mathscr I}/{\mathscr I}^2$ of these left and right ${\mathscr O}$-module structures are equivalent.


\smallskip
\noindent
{\it The isomorphism ${\mathscr O}\oplus \Omega^1\simeq {\mathscr O}_{Y^ {(1)}}$.}
The cotangent sheaf $\Omega^1:=\Omega^1_Y$ can be identified with the pullback $j^{-1}({\mathscr I}/{\mathscr I}^2)$ \cite[p. 407]{Eisenbud}.
This identification is defined by the correspondence
\begin{equation}\label{Omega1Y}
\xi:j^{-1}({\mathscr I}/{\mathscr I}^2)\to \Omega^1,\;\;\;     g+{\mathscr I}^2\mapsto\big(d_xg  \big)|_Y,
\end{equation}
 where $d_x g$ is the exterior derivative of $g$ with respect to the variables $x$; i.e. considering $g$ as function of the variables $x$ and keeping the variables $y$ constant.

\smallskip
 If $l\in\Hat{\mathscr O}$, 
  the correspondence $l+{\mathscr I}^2\mapsto l\circ j\oplus (d_x l)|_Y$
defines 
an isomorphism of right ${\mathscr O}$-modules
\begin{equation}\label{Omega1Yoplus}
j^{-1}(\Hat{\mathscr O} /{\mathscr I}^2)\overset{m}{\simeq} {\mathscr O}\oplus \Omega^1.
\end{equation}

When the abelian sheaf ${\mathscr O}\oplus \Omega^1$ is endowed with the left ${\mathscr O}$-action
\begin{equation}\label{leftprod}
f\cdot(h\oplus \alpha):=fh\oplus (hdf+f\alpha),
\end{equation}
 where $\alpha\in\Omega^1$, then  $m$ in
(\ref{Omega1Yoplus}) is an isomorphism of left ${\mathscr O}$-modules. 
We summarize those well-known results in the following proposition.
\begin{Prop}\label{P:Resumen} With the notations above introduced:
\begin{enumerate}
\item The ${\mathscr O}$-modules $j^{-1}({\mathscr I}/{\mathscr I}^2)$ and $\Omega^1$ are canonically isomorphic.
 \item The correspondence  (\ref{Omega1Yoplus}) defines an isomorphism between  the right  ${\mathcal O}$-modules  $j^{-1}(\Hat{\mathscr O}/{\mathscr I}^2)$ and ${\mathcal O}\oplus{\Omega}^1$. 
\item Equipped ${\mathcal O}\oplus{\Omega}^1$ with the left ${\mathscr O}$-structure defined in (\ref{leftprod}), the correspondence (\ref{Omega1Yoplus})  is also an isomorphism of left ${\mathscr O}$-modules.
\end{enumerate}
\end{Prop}

The exact sequence of right $\Hat{\mathscr O}$-modules
$$0\to {\mathscr I}/{\mathscr I}^2\to\Hat{\mathscr O}/{\mathscr I}^2\to\Hat{\mathscr O}/{\mathscr I}\to 0$$
 gives rise, by means the functor $j^{-1}$, the exact   sequence of right ${\mathscr O}$-modules
 \begin{equation}\label{jexact}
 0\to j^{-1}({\mathscr I}/{\mathscr I}^2)\to j^{-1}(\Hat{\mathscr O}/{\mathscr I}^2)\to j^{-1}( \Hat{\mathscr O}/{\mathscr I})\to 0,
 \end{equation}
or in other terms $0\to\Omega^1\to{\mathscr O}\oplus\Omega^1\to{\mathscr O}\to 0.$


\smallskip
\subsubsection{The first jet sheaf}
 We denote by $\pi_a=p_a\circ i:\Delta \to Y$, for $a=1,2$. One has the following morphism of sheaves rings over $\Delta$  (see Subsection \ref{Ss:Neighborhood})
$$\pi_a^{-1}{\mathscr O}\longrightarrow {\mathscr O}_{Y^{(1)}}= 
{ \Hat{\mathscr O}   } /{ \mathscr I ^2},\;\;\;\; h\mapsto h\circ \pi_a+{\mathscr I}^2.$$
 In particular, one can   consider ${\mathscr O}_{Y^{(1)}}$ as   a right ${\pi_2^{-1}{\mathscr O}}$-module.

 Given ${\mathscr F}$ a left ${\mathscr O}$-module on $Y,$ its inverse image by $\pi_2$ is left ${\mathscr O}_{Y^{(1)}}$-module
  $$\pi_2^*({\mathscr F})={\mathscr O}_{Y^{(1)}}\otimes_{\pi_2^{-1}{\mathscr O}}\pi_2^{-1}{\mathscr F}.$$
And the first jet sheaf ${\mathscr J}^1({\mathscr F})$ of ${\mathscr F}$ is left ${\mathscr O}$-module defined by
\begin{equation}\label{jet1bundle}
{\mathscr J}^1({\mathscr F})={\pi_1}_*\pi_2^*({\mathscr F}).
\end{equation}

Since $\pi_1^*$ is the left adjoint of ${\pi_1}_*$, one has
\begin{equation}\label{adjuntion}
{\rm Hom}\,_{{\mathscr O}}\big({\mathscr F},\,{\mathscr J}^1({\mathscr F})\big)= {\rm Hom}\,_{{\mathscr O}_{Y^{(1)}}} \big(\pi_1^*{\mathscr F}, \,\pi_2^*{\mathscr F}\big).
\end{equation}
 
If ${\mathscr H}$ is a left ${\mathscr O}_{Y^{(1)}}$-module, the ${\mathscr O}$-structure on ${\pi_1}_*{\mathscr H}$ is defined by $f\cdot s=(f\circ\pi_1)s$, where $f\in{\mathscr O}$ and $s\in{\mathscr H}$.
On the other hand, the ${\mathscr O}$ action on 
$$k^*{\mathscr H}={\mathscr O}\otimes_{k^{-1}{\mathscr O}_{\Delta}}k^{-1}{\mathscr H}$$
is determined by 
$f\cdot(1\otimes s)=f\otimes s=1\otimes gs=1\otimes (f\circ\pi_1)s$, where $f=g\circ k$.
Hence, $k^*{\mathscr H}$ and ${\pi_1}_*{\mathscr H}$ are isomorphic ${\mathscr O}$-modules.
 Thus, by (\ref{Omega1Yoplus})
   \begin{equation}\label{J1bis}{\mathscr J}^1({\mathscr F})=k^{-1}({\mathscr O}_{Y^{(1)}})\otimes_{{\mathscr O}}{\mathscr F}=
   j^{-1}(\Hat{\mathscr O}/{\mathscr I}^2)\otimes_{{\mathscr O}}{\mathscr F}.
   \end{equation}

By Proposition \ref{P:Resumen}, ${\mathscr J}^1({\mathscr F})$ is de abelian sheaf ${\mathscr F}\oplus \Omega^1({\mathscr F})$ endowed with the following left ${\mathscr O}$-module structure
\begin{equation}\label{productbyf}
 f(\sigma\oplus\beta)=f\sigma\oplus(f\beta+df\otimes\sigma).
\end{equation}

One has the morphism of abelian sheaves 
 $$\eta:{\mathscr F}\to {\mathscr J}^1({\mathscr F})= {\mathscr F}\oplus\Omega^1({\mathscr F}), 
\;\; \sigma\mapsto \sigma\oplus 0.$$
And from (\ref{productbyf}), it follows 
\begin{equation}\label{etadef}
f\eta(\sigma)= \eta(f\sigma)+df\otimes\sigma.
\end{equation}

Summarizing, one has the exact sequence of ${\mathscr O}$-modules
\begin{equation}\label{exactJ1F}
0\to\Omega^1({\mathscr F})\overset{\iota}\to{\mathscr F}\oplus \Omega^1({\mathscr F})\overset{\pi}{\to} {\mathscr F}\to 0,
\end{equation}
where $\oplus$ is the direct sum in the category of abelian sheaves and
 the left ${\mathscr O}$-action on the central term is defined according to (\ref{productbyf}).

 Thus, taking into account (\ref{Deligne_Gauge}) and (\ref{adjuntion}),
   we give the following definition.
  
  \begin{Def} A holomorphic gauge  field on the coherent sheaf ${\mathscr F}$ is
  an element 
  ${\rm Hom}({\mathscr F},\,{\mathscr J}^1({\mathscr F}))$  that is a right inverse 
    of the ${\mathscr O}$-module morphism  $\pi:{\mathscr J}^1({\mathscr F}) \to {\mathscr F}.$ 
\end{Def}


  \subsubsection{The Atiyah class.}
  The exact sequence (\ref{exactJ1F}), define an element $a({\mathscr F})$ of ${\rm Ext}^1({\mathscr F},\, \Omega^1({\mathscr F}))$, called the Atiyah class of ${\mathscr F}$.
The exact sequence (\ref{exactJ1F}) does not split, in general,  in the category of ${\mathscr O}$-modules. Thus, there will be exist gauge fields on ${\mathscr F}$ iff the Atiyah class $a({\mathscr F})$ vanishes. 

 If there exists   a right inverse $\psi$ of $\pi$; 
then  $\pi(\psi-\eta)=0$; and thus $\nabla:=\psi-\eta$ takes its values in $\Omega^1({\mathscr F})$.
Moreover, by (\ref{etadef}), for $f\in{\mathscr O}$ and $\sigma\in{\mathscr F}$
\begin{equation}\label{LeibnizRule}
 \nabla(f\sigma)=f\psi(\sigma) -f\eta(\sigma)+df\otimes \sigma =f\nabla(\sigma)+df\otimes \sigma;
\end{equation}
that is, for $\nabla$ holds   the Leibniz's rule.

 We denote by ${\mathcal C}$ the set consisting of all  morphisms of  ${\mathbb C}_Y$-modules
$\nabla:{\mathscr F}\to \Omega^1({\mathscr F}) $  satisfying (\ref{LeibnizRule}).
An element $\nabla\in{\mathcal C}$  determines the morphism of ${\mathscr O}$-modules 
$$\kappa: {\mathscr F}\oplus \Omega^1({\mathscr F})\to\Omega^1({\mathscr F}),\;\; (\sigma,\beta)\mapsto \beta-\nabla(\sigma).$$
$\kappa$ is a left inverse of $\iota$. Hence, $\kappa$ defines a splitting of the extension (\ref{exactJ1F}) and thus it determines
a right inverse of $\pi$; i. e. a holomorphic gauge field $\psi$ on ${\mathscr F}$. We have proved the following proposition.
\begin{Prop}\label{P:nabla_psi}
If $a({\mathscr F})=0$, then the map $\nabla\mapsto \nabla+\eta$ defines a bijective correspondence between  ${\mathcal C}$ and the set of holomorphic gauge fields on the coherent sheaf ${\mathscr F}$.
\end{Prop}

\begin{Prop} \label{P:Affine} The set of holomorphic gauge fields on the coherent sheaf ${\mathscr F}$, if it is nonempty,  is an affine space associated to the finite dimensional vector space $\Gamma\big(Y,\,{\mathscr Hom}({\mathscr F},\Omega^1({\mathscr F}))\big).$  
\end{Prop}
{\it Proof.}
The exact sequence (\ref{exactJ1F}) gives rise to corresponding ${\rm Ext}$-sequence
$$0\to{\rm Hom}({\mathscr F},\Omega^1({\mathscr F}))\overset{\lambda}{\to} {\rm Hom}({\mathscr F},{\mathscr J}^1({\mathscr F}))
\overset{\mu}{\to} {\rm Hom}({\mathscr F},{\mathscr F})\overset{\nu}{\to}{\rm Ext}^1({\mathscr F},\Omega^1({\mathscr F})).$$
Since $\mu(\phi)=\pi\circ \phi$, the existence of a holomorphic gauge field $\psi$ on ${\mathscr F}$, is equivalent to
$1_{\mathscr F}\in{\rm im}(\mu)={\ker}(\nu)$. In fact, the Atiyah class $a({\mathscr F})$ is the image of $1_{\mathscr F}$ by $\nu$. 

If $\psi$ and $\psi_1$ are gauge fields on the coherent sheaf ${\mathscr F}$, then $\mu(\psi_1-\psi)=0$; i. e. $\psi_1-\psi\in {\rm im}(\lambda)$. Thus, the set of holomorphic gauge fields on ${\mathscr F}$, if nonempty, is an affine space with vector space  
 ${\rm Hom}({\mathscr F},\Omega^1({\mathscr F}))=\Gamma(Y,\,{\mathscr Hom}({\mathscr F},\Omega^1({\mathscr F})).$ 
  As ${\mathscr Hom}({\mathscr F},\Omega^1({\mathscr F}))$ is a coherent sheaf, the space of holomorphic gauge fields on   ${\mathscr F}$, is a finite dimensional affine space \cite[page 700 ]{G-H}. 
 \qed


  

  \subsubsection{Flat gauge fields}\label{Ss:Flat gauge fields}

Given $\nabla$ a holomorphic connection on $\mathscr F$, it defines a morphism of ${\mathbb C}_Y$-modules $\nabla^{(k)}:\Omega^k({\mathscr F}) {\to} \Omega^{k+1}({\mathscr F})$ in the usual way.  
The composition ${\mathcal K}_{\nabla}:=\nabla^{(1)}\circ\nabla:{\mathscr F}\to\Omega^2({\mathscr F})$ is the curvature of $\nabla$; furthermore,
\begin{equation}\label{mathcalKnabla}
{\mathcal K}_{\nabla}\in{\rm Hom}({\mathscr F},\,\Omega^2({\mathscr F}))=\Gamma(Y,\,{\mathscr Hom}({\mathscr F},\, \Omega^2({\mathscr F}))).
\end{equation}
The connection is said to be {\it flat} if ${\mathcal K}_{\nabla}=0$. In this case, one has the complex
\begin{equation}\label{Complex1}
\Omega^{\bullet}({\mathscr F}):\;\; \Omega^{0}({\mathscr F})\overset{\nabla}{\to} \Omega^{1}({\mathscr F})
\overset{{\nabla}^{(1)}}{\to}
\Omega^{2}({\mathscr F})\to
\end{equation}

\smallskip

A homorphic connection  (not necessarily flat) on a locally free sheaf  ${\mathscr F}$ determines a ${\mathbb C}$-linear map
$$\nabla:{\mathscr A}({\mathscr F})\to  {\mathscr A}^{1,0}({\mathscr F}),$$
satisfying $\nabla(f\tau)=\partial( f)\tau+f\nabla\tau$, for $f$ a smooth function and $\tau$ a section of ${\mathscr A}({\mathscr F})$. In fact, given $\varphi$ is a section of ${\mathscr A}({\mathscr F})$, let $s=\{s_a\}_a$ be a local local frame for
${\mathscr F}$, then $\varphi=\sum_a f^a s_a$, where the $f^a$ are smooth functions. We set
$\nabla(\varphi):=\sum_a(\partial( f^a)s_a +f^a\nabla s_a).$ This is a well-defined section of ${\mathscr A}^{1,0}({\mathscr F})$ (independent of the chosen frame $s$).

 If $\nabla$ is {\it flat}, the extension of $\nabla$ allows us to define the complex
\begin{equation}\label{Complex2}
{\mathscr A}^{\bullet, 0}({\mathscr F}):\;\; {\mathscr A}^{0, 0}({\mathscr F})\overset{\nabla}{\to} {\mathscr A}^{1, 0}({\mathscr F})\overset{\nabla}{\to} {\mathscr A}^{2, 0}({\mathscr F}){\to}
\end{equation}

Analogously, by means of $\nabla+\bar\partial$,  one can construct the following complex, assumed that $\nabla$ is flat,
\begin{equation}\label{Complex3}
{\mathscr A}^{\bullet}({\mathscr F}):\;\; {\mathscr A}^{0}({\mathscr F})\overset{\nabla+\bar\partial}{\longrightarrow} {\mathscr A}^{1}({\mathscr F}) \overset{\nabla+\bar\partial}{\longrightarrow} {\mathscr A}^{1}({\mathscr F})\longrightarrow
\end{equation}

Theorem \ref{IndexTheorem} gives the values of the Euler-Poincar\'e characteristic of the complexes (\ref{Complex1}), (\ref{Complex2}) and (\ref{Complex3}), defined from a holomorphic flat gauge field on the locally free sheaf ${\mathscr F}$.

\medskip

\noindent
{\bf Proof of Theorem 1.} As the dimension is an Euler-Poincar\'e mapping, one has
\begin{equation}\label{Euler-Poin1}
\chi({\mathscr C}^{\bullet})=\sum_i(-1)^i{\rm dim}\,H^i(\Gamma(Y,\,{\mathscr C}^{\bullet})).
\end{equation}

The $p$-column of the following commutative diagram is a fine resolution of  $\Omega^p({\mathscr F}).$
$$\xymatrix{
\Omega^0({\mathscr F}) \ar[d]^i \ar[r]^{\nabla} & \Omega^1({\mathscr F})\ar[d]^i\ar[r]^{\nabla} &  \Omega^2({\mathscr F})\ar[d]^i\ar[r]^{\nabla}   & {}\\
{\mathscr A}({\mathscr F}) \ar[d]^{\bar\partial}\ar[r]^{\nabla} & {\mathscr A}^{1,0}({\mathscr F})\ar[d]^{-\bar\partial}\ar[r]^{\nabla} &  {\mathscr A}^{2,0}({\mathscr F})\ar[d]^{\bar\partial}\ar[r]^{\nabla}    & {} \\
{\mathscr A}^{0,1}({\mathscr F}) \ar[d]^{\bar\partial} \ar[r]^{\nabla} &{\mathscr A}^{1,1}({\mathscr F}) \ar[d]^{-\bar\partial}\ar[r]^{\nabla}   &  {\mathscr A}^{2,1}({\mathscr F}) \ar[d]^{\bar\partial}\ar[r]^{\nabla}     & {} \\
  {} & {} & {}
  } 
  $$
Hence, the total 
complex $({\mathscr A}^{\bullet}({\mathscr F}),\,  D:=\nabla+\bar\partial)$ defined by the double complex $({\mathscr A}^{p,q}({\mathscr F}),\,\nabla, (-1)^p\bar\partial)$ is $q$-isomorphic to the complex in the
top row of diagram, i. e., the complex (\ref{Complex1}). But that total  complex is just (\ref{Complex3}).  Hence,
\begin{equation}\label{H=H}
H^k(\Gamma(Y,\,\Omega^{\bullet}({\mathscr F})))\simeq H^k(\Gamma(Y,\,{\mathscr A}^{\bullet}({\mathscr F}))).
\end{equation}

The complex $\Gamma(Y,\,{\mathscr A}^{\bullet}({\mathscr F}))$, of global sections of (\ref{Complex3})
is elliptic, since the principal symbol of the operator $D=\nabla+\bar\partial$ is equal to the one  of the exterior derivative, $d$. 

The K\"ahler metric and the Hermitian structure determine, in the usual way, an inner product $\diamond$ in the spaces $\Gamma(Y,\,{\mathscr A}^k({\mathscr F}))$. We denote by $D^{\dag}$   the adjoint to $D$ with respect to this inner product.
We can apply the index formula to the operator  operator 
$$P=D+D^{\dag}:\Gamma(Y,\,{\mathscr A}^{\rm even}( {\mathscr F}))\to \Gamma(Y,\,{\mathscr A}^{\rm odd}( {\mathscr F})). $$
  The characteristic classes of ${\mathscr F}$  vanish, since it admits a holomorphic gauge field \cite[Theorem 6]{Atiyah}; thus the Chern character ${\rm ch}({\mathscr F})=r$. Then, by the index formula (see \cite[page 21]{Shanahan})
\begin{equation}\label{EulerCha}
{\rm Index}\,P=r\,{\rm e}(Y).
\end{equation}

On the other hand, the operator $\Delta:=P^{\dagger}P =P P^{\dagger}= D^{\dagger}D +DD^{\dagger},$
$$\Delta\,
: \bigoplus_k\Gamma(Y,\, {\mathscr  A}^k( {\mathscr F}))\to
 \bigoplus_k\Gamma(Y,\, {\mathscr  A}^k( {\mathscr F}))$$
is self-adjoint. Thus,
$$\bigoplus_k\Gamma(Y,\, {\mathscr  A}^k( {\mathscr F}))= 
 {\rm ker}(\Delta)\oplus{\rm im}(\Delta),$$
  and this decomposition is orthogonal \cite[Theorem 5.5, Chap III]{Law-M}. It is easy to check that
 $ {\rm im}(\Delta)$, in turn admits the orthogonal decomposition ${\rm im}( D)\oplus\, {\rm im}(D^{\dagger}).$
 Reasoning as in Hodge decomposition,  it is shown that
\begin{equation}\label{Hodge}
H^i ( \Gamma(Y,\, {\mathscr  A}^{\bullet}( {\mathscr F})))\simeq\{\varphi\in \Gamma(Y,\, {\mathscr  A}^i( {\mathscr F}))\,|\,\Delta\varphi=0\}.
\end{equation}

If $P^{\dagger}P\varphi=0$, then $0=(P^{\dagger}P\varphi)\diamond\varphi=(P\varphi)\diamond (P\varphi)$; so $P\varphi=0$. That is, ${\rm ker}(P^{\dagger}P)={\rm ker}(P).$ Thus,
by (\ref{Hodge})
 $${\rm ker}(P)=
\{\varphi\in\Gamma(Y,\,{\mathscr A}^{\rm even}( {\mathscr F}))\,|\, \Delta \varphi=0\}
\simeq\bigoplus_kH^{2k}(\Gamma(Y,\,{\mathscr A}^{\bullet}({\mathscr F}))).$$

 Analogously, ${\rm coker}(P)={\rm ker}(P^{\dagger})\simeq \oplus_kH^{2k+1}(\Gamma(Y,\,{\mathscr A}^{\bullet}({\mathscr F})))$. 
From  (\ref {Euler-Poin1}) together with (\ref{EulerCha}) and (\ref{H=H}), it follows the first assertion of theorem.

\smallskip

For the case of complex (\ref{Complex2}), we set $Q$ for the operator
$$Q=\nabla+\nabla^{\dag}:\Gamma(Y,\,{\mathscr A}^{{\rm even},0}( {\mathscr F}))\to 
\Gamma(Y,\,{\mathscr A}^{{\rm odd},0}( {\mathscr F})).$$
By de index formula \cite[page 28]{Shanahan} \cite[page 258]{Law-M}
\begin{equation}\label{IndexQ}
{\rm Index}\,Q=(-1)^n\Big(\frac{{\rm ch}\big(\sum_i(-1)^i{\mathscr A}^{i,0}({\mathscr F})\big)\,{\rm td}(T\oplus \bar T )  }{{\rm e}(T)}\Big)[Y].
\end{equation}
The Euler classes of  the holomorphic and antiholomorphic tangent bundles to $Y$ satisfies ${\rm e}(T)=(-1)^n\,{\rm e}(\bar T)$.
Furthermore, ${\mathscr A}^{i,0}\simeq \Lambda^i \bar T$ and (see \cite[page 242]{Law-M})
$${\rm ch}\big( \sum_i(-1)^i\Lambda^i\bar T  \big)=(-1)^n\,{\rm e}(\bar T)\,{\rm td}(T)^{-1}.$$
Thus, (\ref{IndexQ}) reduces  to  ${\rm Index}\,Q=(-1)^n r\,{\rm td}(\bar Y).$

As in the case of the above operator $P$, one has
$${\rm ker}(Q)={\rm ker}(Q^{\dag}Q)=\bigoplus_k H^{2k}(\Gamma(Y,\,{\mathscr A}^{\bullet,0}({\mathscr F}))),$$
and similarly for ${\rm coker}(Q)$. Then the theorem follows.
\qed


\subsection{Holomorphic Yang-Mills fields.}\label{ss:holomorphicYang}

In general, the fiber at $x\in Y$ of a coherent ${\mathscr O}$-module ${\mathscr G}$ 
will be denoted by
${\mathscr G}_{(x)}:= {\mathscr G}_{x}/{\mathfrak m}_x {\mathscr G}_{x}$, where ${\mathfrak m}_x$ is the maximal ideal of ${\mathscr O}_x$.  If ${\mathcal Z}$ is a section of ${\mathscr G}$ the corresponding vector in 
	${\mathscr G}_{(x)}$ is denoted  by $Z(x)$.  
	
	The singularity set ${\mathcal S}$ of ${\mathscr G}$ is an analytic subset of $Y$ whose codimension is greater or equal to $1$. 
Moreover,  ${\mathscr G}$ is locally free on $Y\setminus{\mathcal S}$. 
We set ${ G}$ for the vector bundle over $Y\setminus {\mathcal S}$ with fibers $ G(x):={\mathscr G}_{(x)}$, determined by the locally free sheaf ${\mathscr G}|_{Y\setminus {\mathcal S}}$.

 \begin{Def}\label{D:HermitianSheaf}
A Hermitian  metric on the coherent sheaf ${\mathscr G}$ is a set $\{\langle\,,\,\rangle_x\}_{x\in Y}$ of
Hermitian metrics on the fibers of ${\mathscr G}$, such that, for ${\mathcal Z}_1,{\mathcal Z}_2$ sections of ${\mathscr F}$ on an open $U$ of $Y$,  
the map 
$x\in U\mapsto \langle Z_1(x),\,Z_2(x)\rangle_{x}$
is 
 $C^{\infty}.$   A   sheaf endowed with a Hermitian metric is called   a Hermitian sheaf.
\end{Def}
	
If $Y$ is a K\"ahler manifold and ${\mathcal Z}_1,{\mathcal Z}_2\in\Gamma(Y,\,{\mathscr G})$, we set
\begin{equation}\label{(mathcalZ,Z}
({\mathcal Z}_1,\,{\mathcal Z}_2)=\int_Y\langle Z_1(x),\,Z_2(x)\rangle_{x}\,{\rm d}{\rm vol}=\int_{Y\setminus {\mathcal S}}\langle Z_1(x),\,Z_2(x)\rangle_{x}\,{\rm d}{\rm vol}.
\end{equation}

 Let ${\mathscr F}$ be a coherent sheaf on $Y$.
	For each $x\in Y,$ we denote by $\alpha_x$ and $\lambda_x$ the natural morphisms 
	$$\big({\mathscr Hom}_{\mathscr O}({\mathscr F},\, \Omega^k({\mathscr F}))\big)_x \overset{\alpha_x}{\to}
	{\rm Hom}_{{\mathscr O}_x}\big({\mathscr F}_x,\,\Omega^k_x\otimes_{{\mathscr O}_x} {\mathscr F}_x \big)  \overset{\lambda_x}{\leftarrow}
	\Omega^k_x\otimes_{{\mathscr O}_x}{\rm End}_{{\mathscr O}_x}({\mathscr F}_x).$$
	As ${\mathscr F}$ is coherent, $\alpha_x$ is isomorphism \cite[page 239]{G-R}. Furthermore, if ${\mathscr F}_x $ is free, then $\lambda_x$ is bijective. Hence, for each point $x$ outside of the singularity set ${\mathcal S}$ of ${\mathscr F}$,  the fibre of ${\mathscr Hom}_{\mathscr O}({\mathscr F},\, \Omega^k({\mathscr F})) $ at  $x$
	can be identified with the vector space 
	$\Omega^k{(x)}\otimes{\rm End}({ F}{(x)}).$

	In particular,
	the curvature ${\mathcal K}:={\mathcal K}_{\nabla}$ of a holomorphic connection $\nabla,$ defined in (\ref{mathcalKnabla}), determines the vector
	\begin{equation}\label{tildeKx}
	K(x)\in \Omega^2(x)\otimes{\rm End}( F(x))
	\end{equation}
	for each $x\in Y\setminus{\mathcal S}$. That is, $K$ is a $2$-form ${\rm End}(F)$-valued. 
	
	If ${\mathscr F}$ is a Hermitian sheaf. The K\"ahler structure on $Y$ and metric Hermitian on ${\mathscr F}$   induce a metric on ${\Omega^2\otimes_{\mathscr O}\rm End}(F)$, which will also denoted
$\langle\, ,\,\rangle$. 
According to (\ref{(mathcalZ,Z}),
one defines 
\begin{equation}\label{Norm_Curvature}
\norm{ {\mathcal K}_{\nabla}}^2= \int_{Y\setminus{\mathcal S}}\langle K_{\nabla},\, K_{\nabla}\rangle\,{\rm d}{\rm vol}=
\int_{Y\setminus{\mathcal S}}   |K_{\nabla}\wedge\star\,   K_{\nabla}|,
\end{equation}
where  $|\,\cdot\,|$ is the corresponding norm on ${\rm End}(F)$ and $\star$ is the Hodge star operator.

 More concretely, if locally $  K_{\nabla}$ can be expressed as $\alpha \otimes A$, with $\alpha$ a $2$-form and $A$ a local section of ${\rm End}(F)$, then
the integrand in (\ref{Norm_Curvature}) is $(\alpha\wedge\star\alpha)\langle A\circ A\rangle$. In a local unitary frame of ${\rm End}(F)$, if the connection is compatilble  with the metric,  the matrix $\check A$ associated to $A$ is antihermitian and  
$\langle A\circ A\rangle=- {\rm tr}(\check A\check A)$. That is,
\begin{equation}\label{traza-Norma}
 |  K_{\nabla}\wedge\star\,   K_{\nabla} |=-{\rm tr}\big(  K_{\nabla}\wedge\star\,   K_{\nabla}\big).
\end{equation}

Assumed the set of holomorphic gauge fields on ${\mathscr F}$ is nonempty, the correspondence
$$\nabla\in\{\text{holomorphic gauge fields on}\,{\mathscr F}\}\mapsto \mathcal{ YM}(\nabla)=\norm{{\mathcal K}_{\nabla}}^2$$
is the Yang-Mills' functional \cite[page 417]{Hamilton} \cite[page 44]{Moore} \cite[page 357]{Naber}. The fields on which $\mathcal {YM}$ vanishes are the vacuum states of the corresponding Yang-Mills theory \cite[page 447]{Hamilton}. The gauge fields, where the functional takes a stationary value are the {\it holomorphic Yang-Mills fields}. 

If $\nabla$ is a vacuum state, by (\ref{Norm_Curvature}) it follows $K_{\nabla}=0$, and from the Nakayama's lemma one deduces   ${\mathcal K}_{\nabla}=0$; that is, $\nabla$ is a flat connection.


 \smallskip
 
 Let us assume  that the  sheaf ${\mathscr F}$ admits a holomorphic gauge field $\nabla_0$. By Proposition \ref{P:Affine}, given 
	${\mathcal E}_1,\dots {\mathcal E}_m$, a basis of ${\rm Hom}({\mathscr F},\,\Omega^1({\mathscr F}))$,
	 any holomorphic gauge field can be written
	$\nabla=\nabla_0+\sum\lambda_i{\mathcal E}_i,$ with $\lambda_i\in{\mathbb C}$. The curvature
	$${\mathcal K}_{\nabla}=\nabla\circ\nabla={\mathcal K}_{\nabla_0} +\sum_i\lambda_i {\mathcal B}_i+\sum_{ij}\lambda_i\lambda_j{\mathcal  B}_{ij},$$
	where the ${\mathcal B}$'s are elements of ${\rm Hom}({\mathscr F},\,\Omega^2({\mathscr F}))$
	Thus, 
	$$\norm{{\mathcal K}_{\nabla}  }^2=({\mathcal K}_{\nabla},\,{\mathcal K}_{\nabla})=P(\lambda_1,\dots,\lambda_m),$$
	where $P$ is a polynomial 
 of degree $\leq 4$ in the variables $\lambda_i$. 
	
	\smallskip
	
 \noindent
 {\bf Proof of Theorem 2.}
	The Yang-Mills fields are those $\nabla$ defined by   constants $\lambda_i$ which satisfy the algebraic equations of degree $\leq 3$
 \begin{equation}\label{partialP}
 \frac{\partial\,P}{\partial \lambda_i}=0,\;\;\; i=1,\dots,m.
  \end{equation}

The case $m=2$ is a consequence of B\'ezout's theorem. \qed



\subsubsection{Reflexive sheaves.}

When $Y$ is a Hodge manifold, then it is a smooth projective variety, according to a well-known Kodaira's theorem. By the GAGA correspondence, the coherent analytic sheaves on $Y$ can be identified with the algebraic ones. 
From now on in this Section  \ref{ss:holomorphicYang}, we assume that $Y$ is a Hodge manifold.

 On the other hand, the reflexive sheaves on an algebraic variety might be thought as ``vector bundles with singularities'' \cite[page 121]{Hartshorne1}. The following properties show that these singularities may be  in some cases ``irrelevant".
If ${\mathscr G}$ is a reflexive sheaf  on the algebraic variety $X$, then the codimension of the singularity set ${\rm Sing}$ of ${\mathscr G}$ is greater than $2$
  \cite[Cor. 1.4]{Hartshorne1}. Hence,  the restriction $\Gamma(X,\,{\mathscr G})\to \Gamma(X\setminus {\rm Sing},\,{\mathscr G})$ is an isomorphism \cite[Prop. 1.11]{Hartshorne2}.

    Let us assume that ${\mathscr F}$ is  reflexive sheaf on the Hodge manifold $Y,$ endowed with a Hermitian metric.Then
 ${\mathscr End}({\mathscr F})$ and ${\mathscr Hom}({\mathscr F},\,\Omega^k\otimes_{\mathscr O}{\mathscr F})$ are also a reflexive sheaves \cite[Chapter V, Proposition (4.15)]{Koba}. On the other hand, if ${\mathcal S}$ is the singularity locus of ${\mathscr F} $, one has the isomorphism 
$${\mathscr Hom}({\mathscr F},\,\Omega^k\otimes_{\mathscr O}{\mathscr F})|_{Y\setminus {\mathcal S}}\simeq \big(\Omega^k\otimes_{\mathscr O}{\mathscr End}({\mathscr F})\big)|_{Y\setminus {\mathcal S}}$$
 Thus, we have  the isomorphisms  
	$$
	\Gamma(Y,\,{\mathscr Hom}({\mathscr F},\,\Omega^k\otimes_{\mathscr O}{\mathscr F}))\simeq 
	\Gamma(Y\setminus{\mathcal S},\,\Omega^k\otimes_{\mathscr O}  {\mathscr End}({\mathscr F}))\simeq \Gamma(Y,\,\Omega^k\otimes_{\mathscr O}  {\mathscr End}({\mathscr F})).$$
	Moreover, these {\it finite dimensional} vector spaces are also isomorphic to the space of global sections
	$\Gamma(Y\setminus {\mathcal S},\, \Omega^k\otimes{\rm End}(F))$  of  the vector bundle $\Omega^k\otimes{\rm End}(F)$.
 In particular, 
	${\mathcal K}_{\nabla}$, the curvature of a holomorphic connection $\nabla$   on ${\mathscr F}$, is determined by  the $2$-form ${\rm End}(F)$-valued $K_{\nabla}$ defined     over  $Y\setminus {\mathcal S}$. 
	
	We denote 
	\begin{equation}\label{(p)nabla}
	{}^{(p)}\nabla:\Gamma(Y\setminus{\mathcal S},\,\Omega^p\otimes_{\mathscr O}{\rm End}({ F}))\to
	\Gamma(Y\setminus{\mathcal S},\,\Omega^{p+1}\otimes_{\mathscr O}{\rm End}({ F})),
	\end{equation}
	the operator defined by the connection $\nabla$. In this notation Bianchi's identity is read as  
	$$ {}^{(2)}\nabla { K}_{\nabla}=0.$$ 


 \smallskip
 If $\nabla$ is a holomorphic gauge field on ${\mathscr F}$, according to Proposition \ref{P:Affine}, any other field is of the form 
$\nabla+{\mathcal E}$, with
 ${\mathcal E}\in{\rm Hom}({\mathscr F},\, \Omega^1\otimes_{\mathscr O} {\mathscr F})$. By the above identifications ${\mathscr E}$ is determined by the corresponding section $E\in \Gamma(Y\setminus {\mathcal S},\, \Omega^1\otimes{\rm End}(F))$.

Considering a ``variation'' $\nabla_{\epsilon}=\nabla+\epsilon{\mathscr E}$ of $\nabla$, then $K_{\nabla_{\epsilon}}=K_{\nabla}+\epsilon\nabla E+O(\epsilon ^2)$, where
 $\nabla   E$ is the covariant derivative of $  E$.

\begin{equation}\label{(1/2)}(1/2)\frac{d}{d\epsilon}\Big|_{\epsilon=0}||{\mathcal K}_{\nabla_{\epsilon}}||^2=\int_{Y\setminus{\mathcal S}}\langle K_{\nabla},\, \nabla  E \rangle\, {\rm d}{\rm vol}=: (  K_{\nabla},\, \nabla{ E}).
\end{equation}
Therefore, $\nabla$ is a Yang-Mills field if for any ``variation'' ${\mathcal E}$ of $\nabla$ 
\begin{equation}\label{YanMillsEq}
(  K_{\nabla},\, \nabla{ E})=0.
\end{equation}
 In particular, the flat holomorphic gauge fields are Yang-Mills.

	
	On the other hand,   the orthogonality condition (\ref{YanMillsEq}) which satisfy the Yang-Mills fields gives rise to the following proposition.
	\begin{Prop}\label{P:ortho}
	The   holomorphic gauge field $\nabla$ on the reflexive sheaf ${\mathscr F}$ is a Yang-Mills field iff
	its curvature $K_{\nabla}\in \Gamma(Y\setminus{\mathcal S},\,\Omega^2\otimes{\rm End}({ F}))$ is orthogonal to the vector space ${\rm im}\,(^{(1)}\nabla)$.
	\end{Prop}
	
	\begin{Cor}\label{Cor:Flat}
	If $H^2(\Gamma(Y,\,\Omega^{\bullet}\otimes{\rm End }(F))=0$, then any Yang-Mills field on ${\mathscr F}$ is flat. 
	\end{Cor}
	{\it Proof.}
	  By Bianchi's identity $K_{\nabla}\in {\rm ker}\,^{(2)}\nabla$. By the hypothesis $K_{\nabla}\in{\rm im}\,^{(1)}\nabla $. If $\nabla$ is a Yang-Mills field, then $K_{\nabla}$
is a vector orthogonal to ${\rm im}\,( ^{(1)}\nabla),$ according to Proposition \ref{P:ortho}. Thus, $K_{\nabla}=0$,  and by Nakayama lemma ${\mathcal K}_{\nabla}=0$.
\qed

A flat connection on a vector  bundle 
defines 
a $D$-module structure
 on the corresponding ${\mathscr O}$-module. Thus, by the Riemann-Hilbert correspondence, one has the following corollary.
\begin{Cor}\label{C:representation} Let ${\mathscr F}$ be a locally free sheaf such that $H^2(\Gamma(Y,\,\Omega^{\bullet}\otimes{\rm End }(F))=0$. If ${\mathscr F}$   admits a holomorphic Yang-Mills field, then it
  is defined by a 
  representation of $\pi_1(Y)$.
\end{Cor}

  Under hypotheses very different from ours, other authors have shown that vector bundles that support holomorphic connections  
  are actually flat vector bundles (see \cite{Biswas, B-I}). 
	
		\smallskip
		
The orthogonality condition (\ref{YanMillsEq}) implies $^{(1)}\nabla^{\dagger}  K_{\nabla}=0$, where
	$$^{(1)}\nabla^{\dagger}: \Gamma(Y\setminus{\mathcal S},\,\Omega^2\otimes_{\mathscr O}{\rm End}({ F}))\to 
	\Gamma(Y\setminus{\mathcal S},\,\Omega^1\otimes_{\mathscr O}{\rm End}({ F}))$$
	is the adjoint of $^{(1)}\nabla$. 
	By the Bianchi's identity,  if $\nabla$ is a Yang-Mills field, then
	\begin{equation}\label{DeltaKnabla}
	^{(2)}\nabla K_{\nabla}=0,\;\; ^{(1)}\nabla^{\dagger}  K_{\nabla}=0,
	\end{equation}
	and conversely.

	
	\subsubsection{The case   ${\rm rank}\,{\mathscr F}=1$.}
	Let us  assume that ${\mathscr F}$ is a locally free sheaf of rank $1$.
	It is known that a necessary and sufficient condition for ${\mathscr F}$ to admit a connection is that F is flat (see Proposition \ref{a=c-1} in Appendix)
	 
	On the other hand, ${\mathscr End}({\mathscr F})$ is the sheaf associated to the trivial line bundle ${\mathbb C}\times Y\to Y.$ Let $s$ be a local frame of the corresponding line bundle $F.$  A  
	holomorphic connection $\nabla$ on ${\mathscr F}$  
	is locally determined by a 
	holomorphic $1$-form $A$, $\nabla s=A s$. In this frame 
	$\nabla(\beta)=\partial\beta+A\wedge\beta-(-1)^p\beta\wedge A=\partial\beta,$
	for any ${\rm End}(F)$-valued $p$-form $\beta$. Thus, the operator (\ref{(p)nabla}) reduces to 
	\begin{equation}\label{pnabla=partial}
	^{(p)}\nabla=\partial:\Gamma(Y,\,\Omega^p)\to\Gamma(Y,\,\Omega^{p+1}).
	\end{equation}
	 The curvature of this connection is given by the holomorphic $2$-form $\partial A$, and   the Bianchi's identity  reduces  to the obvious relation $\partial K_{\nabla}=0$.
	 
	In this case, the complex $(\Omega^{\bullet}\otimes {\rm End}(F),\,\nabla)$ is  the holomorphic de Rham complex 
	$(\Omega^{\bullet},\,\partial)$ \cite[Sect 8.2.1]{Voisin}. Denoting by ${\mathscr A}^{p,q}$ the sheaf of the smooth $(p,q)$-forms on $Y$, then $({\mathscr A}^{p,\bullet},\,\partial)$ is a fine resolution of $\Omega^p$. Thus, as in the proof of Theorem \ref{IndexTheorem}, the total complex
	 associated to the double complex complex $({\mathscr A}^{p,q};\,\partial,\,(-1)^p\bar\partial)$ is quasi-isomorphic to the holomorphic de Rham complex. Since this total complex is precisely the usual de Rham complex, one has
	\begin{equation}\label{HjOmaga}
	H^j(\Gamma(Y,\,\Omega^{\bullet}))=H^j(Y,\,{\mathbb C}).
	\end{equation}


	\smallskip
	
	\noindent	
	{\bf  Proof of Theorem \ref{impatial=0}.}
	 We can assume that it is a line bundle
	\cite[Prop. 1.9]{Hartshorne1}; thus, the singularity locus $\mathcal S$ is the empty set. Since $c_1({\mathscr F})=0$, the Atiyah class vanishes (see Proposition \ref{a=c-1} in Appendix) and ${\mathscr F}$ admits holomorphic connections.

	Let $\nabla$ be a holomorphic gauge field on ${\mathscr F}$. By the Bianchi's identity, 
	$K_{\nabla}$ defines a cohomology class in the space (\ref{HjOmaga}) with $j=2$. 
	Any other gauge field $\tilde\nabla$ is an element
of $\nabla+{\rm Hom}({\mathscr F},\,\Omega^1\otimes_{\mathscr O}{\mathscr F})$.
	As ${\mathscr F}$ is a locally free sheaf with rank $1$
\begin{equation}\label{HomF,Omega1}
{\rm Hom}({\mathscr F},\,\Omega^1\otimes_{\mathscr O}{\mathscr F})\simeq\Gamma\big(Y,\,\Omega^1\otimes_{\mathscr O}{\mathscr End}({\mathscr F})\big) 
\simeq\Gamma(Y,\,\Omega^1).
\end{equation}

	By (\ref{pnabla=partial}), the curvature of 
	$\tilde\nabla$ has the form $K_{\tilde\nabla}=K_{\nabla}+\partial E$, with $E\in \Gamma(Y,\,\Omega^1)$ a holomorphic $1$-form.
	As $\bar\partial E=0$,  the curvatures $K_{\nabla}$ and $K_{\tilde\nabla}$ determine the same cohomology class. We will denote by $\mathbf{c}$ this cohomology class, defined by curvature of any holomorphic gauge field on ${\mathscr F}$.

	If $\widehat\nabla$ is an arbitrary Yang-Mills field, then $K_{\widehat\nabla}$
	satisfies (\ref{DeltaKnabla}); that is, $K_{\widehat\nabla}$ is $\partial$-harmonic. As $Y$ is a K\"ahler manifold,  $K_{\widehat\nabla}$ is also $d$-harmonic. Hence, the norm of $K_{\widehat\nabla}$ minimizes the corresponding norm  in its cohomology class. That is,
	\begin{equation}\label{minimum}
	\mathcal{YM}(\widehat\nabla)=\norm{K_{\widehat\nabla}}^2={\min}\{\norm{\beta}^2\,|\, \beta\in{\bf c}\}.
	\end{equation}
	
	On the other hand, if $\nabla_0$ is a holomorphic gauge field,  $\frac{i}{2\pi}[K_{\nabla_0}]=c_1({\mathscr F})=0.$ That is,
	$$K_{\nabla_0}=(\partial+\bar\partial)(B^{1,0} +B^{0,1} ).$$    
			Since $K_{\nabla_0}$ is a $(2,0)$-form, it follows $\bar\partial B^{1,0}=0$ and $d B^{0,1}=0$. That is, $K_{\nabla_0}=\partial B^{1,0}$.
			The holomorphic connection  $\tilde\nabla:=\nabla_0-B^{1,0}$ has  curvature zero, hence it satisfies (\ref{YanMillsEq}); that is $\tilde\nabla$ is a Yang-Mills field. Therefore the cohomology class ${\mathbf c}=0$. It follows from  (\ref{minimum}), that $\norm{K_{\nabla}}=0$, for any holomorphic Yang-Mills field.
	\qed

\smallskip

\begin{Prop}\label{Divisor}
Let ${\mathscr L}$ be a line bundle with $c_1({\mathscr L})=0$.
If
the Hodge number $h^{1,0}(Y)=1$, then either  the cardinal $\#{\sf YM}({\mathscr L})=1$,  
or  any gauge field on ${\mathscr L}$ is Yang-Mills. The latter case occurs when $H^0(Y,\,\Omega^1)= {\mathbb C}.$
\end{Prop}
{\it Proof.}
The vanishing of the Chern class implies that there exist homolorphic gauge fields on ${\mathscr L}$. Let $\nabla_0$ denote  a holomorphic gauge field. 
 According to (\ref{HomF,Omega1}), any other gauge field is an element
of $\nabla_0+ \Gamma(Y,\,\Omega^1)$

  Since $h^{1,0}(Y)=1$, any holomorphic gauge field on ${\mathscr L}$ is of the form
$\nabla=\nabla_0+\lambda E$, with $\lambda\in{\mathbb C} $ and $0\ne E\in H^0(Y,\,\Omega^1).$ The curvature 
$$K_{\nabla}=K_{\nabla_0}+\lambda\nabla_0(E).$$
Hence the polynomial
$$P(\lambda)\equiv\norm{K_{\nabla} }^2 = \norm{K_{\nabla_0} }^2+2\lambda(K_{\nabla_0},\,\nabla_0(E))+\lambda^2\norm{\nabla_0(E)}^2.$$
If $\nabla_0(E)\ne 0$ the equation $\frac{d\,P}{d\lambda}=0$ has only one solution, so $\#{\sf YM}({\mathscr L})=1$. 

 By contrast, when
 $\nabla_0(E)=0$,  it follows that $P(\lambda)=\norm{K_{\nabla_0} }^2$, for all $\lambda$. Hence, in this case, the Yang-Mills functional is constant; thus, every holomorphic  gauge field is a Yang-Mills field.
 On the other hand, since $0=\nabla_0(E)=\partial E$, then $E$ is constant.
  \qed


\section{Fields on a brane}\label{S:Fields_brane}


\subsection{Gauge fields on a $B$-brane.} \label{Ss:Fields_brane}  Let $({\mathscr F}^{\bullet},\,\delta^{\bullet})$ be a $B$-brane  on the complex manifold $Y$; that is, ${\mathscr F}^{\bullet}$ is an object of the   category $D^b(Y)$, the bounded derived category of coherent sheaves over $Y$. 
According to the observation at the beginning of Subsection \ref{ss:FirstJet} (see (\ref{Deligne_Gauge})), we define a gauge field on the brane ${\mathscr F}^{\bullet}$ as an element of 
${\rm Hom}_{D^b(Y^{(1)})}(L\pi_1^*{\mathscr F}^{\bullet},\, L\pi_2^*{\mathscr F}^{\bullet})$
which lifts the  automorphism identity of $({\mathscr F}^{\bullet},\,\delta^{\bullet})$.
 This condition will be explained below in the rigorous definition of this concept. 

 As $\pi_i$ is flat, $L\pi_i^*$ is the usual inverse image $\pi_i^*$.
 By the adjunction relation
\begin{equation}\label{HomDbY}
{\rm Hom}_{D^b(Y^{(1)})} (\pi_1^*{\mathscr F}^{\bullet},\, \pi_2^*{\mathscr F}^{\bullet} )\simeq 
{\rm Hom}_{D^b(Y)} ( {\mathscr F}^{\bullet},\,{\mathscr J}^1({\mathscr F}^{\bullet})),
\end{equation}    
 where 
$${\mathscr J}^1({\mathscr F}^{\bullet}):=R\pi_{1*} \pi_2^*{\mathscr F}^{\bullet}\simeq{\mathscr O}_{Y^{(1)}}\otimes^L{\mathscr F}^{\bullet}.$$ 

Therefore, the gauge fields on ${\mathscr F}^{\bullet}$ are elements of the group 
$${\rm Ext}^0({\mathscr F}^{\bullet},\,{\mathscr J}^1({\mathscr F}^{\bullet}));$$
i.e., open strings between ${\mathscr F}^{\bullet}$ and ${\mathscr J}^1({\mathscr F}^{\bullet})$ with ghost
 number $0$ \cite[Sect. 5.2]{Aspin}, \cite{Katz-Sharpe}.

As ${\mathscr O}_{Y^{(1)}}$ is the locally free module ${\mathscr O}\oplus \Omega^1$, then
 ${\mathscr J}^1({\mathscr F}^{\bullet})\simeq{\mathscr F}^{\bullet}\oplus \Omega^1({\mathscr F}^{\bullet)})$, with ${\mathscr O}$-structure given by (see (\ref{productbyf}))
\begin{equation}\label{fsigmabullet}
f\ast(\sigma^{\bullet}\oplus \beta^{\bullet})=f\sigma^{\bullet}\oplus(f\beta^{\bullet}+df\otimes\sigma^{\bullet}).
\end{equation}

The exact sequence of complexes of ${\mathscr O}$-modules 
$$0\to\Omega^1({\mathscr F}^{\bullet})\overset{\iota}{\rightarrow}{\mathscr J}^1({\mathscr F}^{\bullet} ) \overset{\pi}{\rightarrow} {\mathscr F}^{\bullet}\to 0$$
determines a distinguished triangle
$$\Omega^1({\mathscr F}^{\bullet})\overset{\iota}{\rightarrow}{\mathscr J}^1({\mathscr F}^{\bullet} ) \overset{\pi}{\rightarrow} {\mathscr F}^{\bullet}\overset{+1}{\to}$$
in the category $D^b(Y)$ \cite[page 46]{Kas-Sch}, \cite[page 157]{Ge-Ma}. As   ${\rm Hom}_{D^b(Y)}({\mathscr F}^{\bullet},\,.\,)$ is a cohomological functor, it follows
that 
\begin{align}&{\rm Hom}_{D^b(Y)}({\mathscr F}^{\bullet},\, \Omega^1({\mathscr F}^{\bullet}))\to
 {\rm Hom}_{D^b(Y)}({\mathscr F}^{\bullet},\, {\mathscr J}^1({\mathscr F}^{\bullet}))\overset{\pi\circ}{\longrightarrow}
{\rm Hom}_{D^b(Y)}({\mathscr F}^{\bullet},\,{\mathscr F}^{\bullet}) \notag \\ 
&\to {\rm Ext}^1({\mathscr F}^{\bullet},\, \Omega^1({\mathscr F}^{\bullet}))\to  \notag
\end{align}
is an exact sequence. The Atiyah class of ${\mathscr F}^{\bullet}$  is the image of $1\in {\rm Hom}_{D^b(Y)}({\mathscr F}^{\bullet},\,{\mathscr F}^{\bullet})$  in ${\rm Ext}^1({\mathscr F}^{\bullet},\, \Omega^1({\mathscr F}^{\bullet}))$. Thus, we give the following definition.

 \begin{Def}\label{D:gauge} A gauge field on ${\mathscr F}^{\bullet}$ is an element $\psi\in{\rm Hom}_{D^b(Y)}({\mathscr F}^{\bullet},\, {\mathscr J}^1({\mathscr F}^{\bullet}))$, such that $ {\pi\circ}\psi =1\in {\rm Hom}_{D^b(Y)}({\mathscr F}^{\bullet},\,{\mathscr F}^{\bullet})$ 
\end{Def} 
 From the above exact sequence, it follows the following proposition.
\begin{Prop}\label{P:Atiyah}
The vanishing of the Atiyah class ${\mathscr F}^{\bullet}$ is a necessary and sufficient condition for the existence of gauge fields on this brane.  Furthermore, the set of gauge fields on  ${\mathscr F}^{\bullet}$, if is nonempty, is an affine space over the finite dimensional vector space ${\rm Ext}^0({\mathscr F}^{\bullet},\,\Omega^1({\mathscr F}^{\bullet}))$.
\end{Prop}


 \subsubsection{Gauge fields on $B$-branes over ${\mathbb P}^n$}\label{Ss: CPn}
According to the Beilinson spectral sequence \cite[Chap. 2, Sect. 3.1]{O-S-S} 
the set (\ref{set})
is a strong complete exceptional sequence in the derived category  $D^b({\mathbb P}^n)$ \cite[Sect 8.3]{Huybrechts}. Thus, $D^b({\mathbb P}^n)$ is equivalent to the smallest triangulated subcategory that contains this exceptional family.

As usual, we denote with ${\mathscr F}^{\bullet}[l]$, with $l\in{\mathbb Z}$, the complex ${\mathscr F}^{\bullet}$ shifted $l$ to the left. 
   Given  ${\mathscr A}, {\mathscr B}$ elements of the generating set (\ref{set}),
  let us consider morphisms $h$ between ${\mathscr A}':={\mathscr A}[l]$ and
  ${\mathscr B}':={\mathscr B}[l']$. We denote by ${\rm Cone}(h)={\mathscr A}'[1]\oplus {\mathscr B}'$ the mapping cone of $h$
   \cite[page 154]{Ge-Ma}. We define ${\sf E}^{(1)}$ the set obtained adding to ${\sf E}$ the elements of the form ${\rm Cone}(h)$. Hence, an element of ${\sf E}^{(1)}$ is a complex whose term at a position $p$ is either  $0$, or ${\mathscr O}(k)$, or a direct sum of ${\mathscr O}(k_1)\oplus {\mathscr O}(k_2)$, with $-n\leq k,k_1,k_2\leq 0$.
   
    Repeating the process with the elements of ${\sf E}^{(1)}$, one obtains ${\sf E}^{(2)}$, etc. The objects of the triangulated subcategory generated by the family ${\sf E}$ are elements which belong to some ${\sf E}^{(m)}$.
      Therefore, an object  of the triangulated subcategory of $D^b({\mathbb P}^n)$  generated by 
			(\ref{set}) is a complex $({\mathscr G}^{\bullet},\,d^{\bullet})$, where ${\mathscr G}^{p}$ is a sheaf  of the form 
 \begin{equation}\label{bigoplus}
 {\mathscr G}^{p}=\bigoplus_{i\in S_p}{\mathscr O}(k_{pi}),
 \end{equation}
 with $-n\leq k_{pi}\leq 0$ and $i$ varying in a finite set $S_p$. (When $i$ ``runs over the empty set'', the direct sum is taken to be $0$). 

\smallskip
\noindent
{\bf Proof of Theorem \ref{Th:card}.}
In general, given two bounded below  complexes $A^{\bullet}$ and  $B^{\bullet}$ in an abelian category ${\mathfrak A}$, the complex
${\rm Hom}^{\bullet}(A^{\bullet},\,B^{\bullet})$ is defined by (see \cite[page 17]{Iversen})
\begin{equation}\label{deltamg}
{\rm Hom}^m(A^{\bullet},\,B^{\bullet})=\prod_{p\in{\mathbb Z}}{\rm Hom}_{\mathfrak A}\big(A^p,\, B^{p+m}   \big),
\end{equation}
with the differential  $\delta_H$.
\begin{equation}\label{deltamg1}
(\delta_H^mg)^p=\delta_B^{m+p}g^p+(-1)^{m+1}g^{p+1}\delta_A^p.
\end{equation}

\smallskip

  As the complex  ${\mathscr G}^{\bullet}$   defined in (\ref{bigoplus})  consists of 
  locally free ${\mathscr O}$-modules, then \cite[Chap III, 6.5.1]{Hartshorne}
 \begin{align} &{\rm Hom}_{D^b({\mathbb P}^n)}\big({\mathscr G}^{\bullet},\,\Omega^1({\mathscr G}^{\bullet})\big)=
 H^0{\rm Hom}^{\bullet}\big({\mathscr G}^{\bullet},\,\Omega^1({\mathscr G}^{\bullet})\big)= \notag \\ \notag
 &\{g\in{\rm Hom}^{0}\big({\mathscr G}^{\bullet},\,\Omega^1({\mathscr G}^{\bullet})\big)\,|\,\delta_H g=0\},
\end{align}
where $\delta_H$ is the operator defined in (\ref{deltamg1}).
 Hence, according to (\ref{deltamg}), it follows
 $${\rm Hom}_{D^b({\mathbb P}^n)}\big({\mathscr G}^{\bullet},\,\Omega^1({\mathscr G}^{\bullet})\big)\subset
  {\rm Hom}^{0}\big({\mathscr G}^{\bullet},\,\Omega^1({\mathscr G}^{\bullet})\big)=
  \prod_p{\rm Hom}({\mathscr G}^p,\,\Omega^1({\mathscr G}^p)).$$
   By the additivity of the functor ${\rm Hom}(\,.\,,\,.\,)$, it follows
  \begin{equation}\label{proplus}
  {\rm Hom}_{D^b({\mathbb P}^n)}\big({\mathscr G}^{\bullet},\,\Omega^1({\mathscr G}^{\bullet})\big)\subset
  \prod_p\bigoplus_{i,j}{\rm Hom}\big( {\mathscr O}(k_{pi}),\,\Omega^1(k_{pj})  \big), \notag
   \end{equation}
	where $\Omega^1(k)$ is the twisted sheaf $\Omega^1\otimes_{\mathscr O}{\mathscr O}(k)$.

 The summand ${\rm Hom}\big( {\mathscr O}(k_{pi}),\,\Omega^1(k_{pj})  \big)$ is equal to
 $${\rm Hom}\big( {\mathscr O},\,\Omega^1(k_{pj}-k_{pi})  \big)=\Gamma({\mathbb P}^n,\,\Omega^1(k_{pj}-k_{pi}))=0,$$
since $H^0({\mathbb P}^n,\,\Omega^1(k))=0$, for any  $k$ \cite[page 4]{O-S-S}. Therefore,
  $${\rm Hom}_{D^b({\mathbb P}^n)}\big({\mathscr G}^{\bullet},\,\Omega^1{\mathscr G}^{\bullet}\big)=0.$$  
From Proposition \ref{P:Atiyah}, it follows   Theorem \ref{Th:card}. \qed

\medskip

\noindent
{\bf Proof of Theorem \ref{Th:2}.}
 Let $\psi$  be a holomorphic gauge field on the above $B$-brane ${\mathscr G}^{\bullet}$. Then
  \begin{align}
  &\psi\in{\rm Hom}_{D^b({\mathbb P}^n) }\big( {\mathscr G}^{\bullet},\,{\mathscr J}({\mathscr G}^{\bullet} ) \big)=
  H^0{\rm Hom}^{\bullet} \big({\mathscr G}^{\bullet},\, {\mathscr J}({\mathscr G}^{\bullet})  \big)\subset \notag \\  
 &\prod_{p}{\rm Hom} \big( {\mathscr G}^p,\, (\Omega^1({\mathscr G}^p) \oplus {\mathscr G}^p) \big). \notag 
 \end{align}
Thus, $\psi$ determies a family $\{\psi^p:  {\mathscr G}^p \to \Omega^1({\mathscr G}^p) \oplus {\mathscr G}^p\} $ of morphisms of ${\mathscr O}$-modules.
As $\psi$ is a right inverse of $\pi$ (Definition \ref{D:gauge}), $ \psi^p=\nabla^p\oplus {\rm id}^p$, where 
$\nabla^p:{\mathscr G}^p \to \Omega^1({\mathscr G}^p)$. 
 The property $\psi(f\sigma)=f\ast\psi(\sigma)$   
implies that $\nabla(f\sigma)=df\sigma+f\nabla(\sigma)$, for $f\in{\mathscr O}$. Hence, 
  $\nabla^p$ is a holomorphic connection on ${\mathscr G}^p$, for any $p$.

 The trace
  of the curvature  of $\nabla^p$ is a holomorphic $2$-form on ${\mathbb P}^n$; as $H^0({\mathbb P}^n,\,\Omega^2)=0$, that trace vanishes. Hence, the first Chern class of vector bundle associated to the locally free  sheaf ${\mathscr G}^p$ vanishes.
  
   On the other hand, the first Chern class of ${\mathscr G}^p$ is the sum  
   $$\sum_{i} c_1({\mathscr O}(k_{pi})).$$
    This class is 0 iff $k_{pi}=0$ for all $i$, since the $k_{pi}\leq 0$. 
 
 Therefore, the existence of a holomorphic gauge field on the brane  ${\mathscr G}^{\bullet}$ defined in (\ref{bigoplus}) implies that 
 ${\mathscr G}^{\bullet}$ is a sequence of direct sum of copies of ${\mathscr O}$
 \begin{equation}\label{bigoplusimp}
 \dots \to \bigoplus_{i\in S_p}{\mathscr O}\overset{d^p}{\longrightarrow} \bigoplus_{i\in S_{p+1}}{\mathscr O}\to\dots
 \end{equation}
   Since ${\rm Hom}_{\mathscr O}({\mathscr O},\, {\mathscr O})\simeq {\mathbb C}$, the map $d^p$ is given by a constant complex matrix. 
   
 On the other hand, given the brane (\ref{bigoplusimp}),  if the set of indices $S_p$ has $m_p$ elements, on  $\oplus_{1}^{m_p}{\mathscr O}$ we define the map $\varphi^p$,   
   $$\varphi^p(\sigma_1\oplus\dots\oplus\sigma_{m_p})= \partial\sigma_1\oplus\dots\oplus\partial\sigma_{m_p}.$$ 
   It is a holomorphic connection on $\oplus_{1}^{m_p}{\mathscr O}$. Moreover, the family $\{\varphi^p\}$ is compatible with the ``constant'' differentials $d^p$. Thus, this family is  a holomorphic gauge field  on the brane defined by the complex
   (\ref{bigoplusimp}).  \qed
   



\subsubsection{\it The homotopy category}\label{homotopycategory}
Let $X$ be a smooth projective variety.
 The category $\mathfrak{Coh}(X)$ of coherent sheaves over $X$  has not enough injectives, for this reason it is convenient to regard $D^b(X)$
 as a subcategory of $D^b({\mathscr O_X})$, the bounded derived category of the ${\mathscr O}_X$-modules.  
In fact, $D^b(X)$ is equivalent to $D^b({\mathscr O}_X)_{\rm coh}$, the full subcategory of $D^b({\mathscr O}_X)$ consisting of the complexes with   coherent cohomology  \cite[Chapter II]{Hartshorne_R} \cite[Exp II]{B-G-I}.
Thus,  $D^b(X)$ can be identified with the homotopy category whose objects are the complexes ${\mathscr G}^{\bullet}$ of injective ${\mathscr O}_X$-modules, such that its cohomogy is bounded and coherent; 
i.e. ${\mathscr H}^j({\mathscr G}^{\bullet})$ is coherent and vanishes for $|j|>>0$. We will denote this homotopy category by $K^b(X)_{\rm coh}$.

\smallskip
Henceforth, we assume that $Y$ is a {\it Hodge manifold}; in this way we can identify coherent analytic sheaves on $Y$ with algebraic ones. 
 	Hence, the $B$-branes on 
	$Y$ can be considered as objects of     $K^b(Y)_{\rm coh}$.  
	 
  	 In accordance with the preceding paragraph, one can assume that the brane $({\mathscr F}^{\bullet},\delta^{\bullet} )$ is a complex of {\it injective} ${\mathscr O}$-modules with coherent cohomology modules satisfying ${\mathscr H}^i({\mathscr F}^{\bullet} )=0$ for $|i|>>0$. Moreover,
	$${\rm Hom}_{D^b(Y^{(1)})} (\pi_1^*{\mathscr F}^{\bullet},\, \pi_2^*{\mathscr F}^{\bullet} )\simeq 
{\rm Hom}_{K^b(Y)_{\rm coh}}({\mathscr F}^{\bullet},\, \widehat{{\mathscr J}^1({\mathscr F}^{\bullet}})),$$ 
where $\widehat{{\mathscr J}^1({\mathscr F}^{\bullet}})$ is an object of $K^b(Y)_{\rm coh}$ $q$-isomorphic to 
${\mathscr J}^1({\mathscr F}^{\bullet})$.

Thus, the elements of the space (\ref{HomDbY}) can be identified with
	morphisms in the homotopyy category $K^b(Y)_{\rm coh}$.
	
	From now on in this subsection,
	 \emph{we delete the bullets in the notation for the complexes} and set ${\mathscr J}:={\mathscr J}^1({\mathscr F}^{\bullet})$.
	Then there is a morphism $\hat\pi : \widehat{{\mathscr J} } \to{\mathscr F}$ in $K^b(Y)_{\rm coh}$ determined by $\pi$,
	and the gauge field $\psi$ can be regarded as a morphism in $K^b(Y)_{\rm coh}$,
$\hat\psi:{\mathscr F}\to\widehat{{\mathscr J}}$, such that $\hat\pi\circ\hat\psi=1$.

On the other hand, $\hat\psi$ as a morphism of a homotopy category, determines a well-defined morphism of ${\mathscr O}$-modules between the cohomologies, that will be denoted in bold,
$$\pmb\psi^j:{\mathscr H}^j({\mathscr F})\to  {\mathscr H}^j(\hat{\mathscr J}) = {\mathscr H}^j({\mathscr J}).$$
 Similarly, one has the canonical projection
 $$\pmb\pi^j: {\mathscr H}^j({\mathscr J}) = {\mathscr H}^j(  {\mathscr F})\oplus\big({\Omega^1}\otimes_{\mathscr O} {\mathscr H}^j(  {\mathscr F})\big)\to {\mathscr H}^j({\mathscr F}),$$
  satisfying $\pmb\pi^j\pmb\psi^j=1$.
 
 We set $\pmb \eta^j$ for the morphism of abelian sheaves defined by the inclusion in the direct sum
 $$\pmb\eta^j:{\mathscr H}^j({\mathscr F}) \rightarrow {\mathscr H}^j( {\mathscr J})={\mathscr H}^j(  {\mathscr F})\oplus\big({\Omega^1}
 \otimes_{\mathscr O} {\mathscr H}^j(  {\mathscr F})\big),$$


  Hence, $\pmb\pi^j(\pmb\psi^j-\pmb\eta^j)=0$, 
  and  thus $\pmb\psi^j-\pmb\eta^j$ defines a morphism of abelian sheaves 
  \begin{equation}\label{vartethaConnections}
	\vartheta^j:{\mathscr H}^j({\mathscr F})\to{\Omega^1}\otimes_{\mathscr O} {\mathscr H}^j({\mathscr F}),
	\end{equation}
which, by (\ref{fsigmabullet}), satisfies  the Leibniz's. That is, 
\begin{Prop}\label{P:Connvartheta} The gauge field $\psi$ on the brane ${\mathscr F}$ determines on each sheaf 
 ${\mathscr H}^j({\mathscr F})$  a connection $\vartheta^j$.
\end{Prop}

\begin{Def}\label{D:flat}
The gauge field $\psi$ is called flat, if the curvature of the connection $\vartheta^j$ vanishes, for all $j$.
\end{Def}

\begin{Rems}\label{R:variation}
Let $\psi,\phi$ be two gauge fields on ${\mathscr F}$, we set 
$$\xi:=\phi-\psi\in{\rm Hom}_{D^b(Y)}\big({\mathscr F},\,{\mathscr J}({\mathscr F})\big),$$
and denote by $\hat\xi$ the corresponding morphism ${\mathscr F}\to \widehat{{\mathscr J}({\mathscr F})}$ in the category $K^b(Y)_{\rm coh}.$ Thus, $\hat\xi$ determines a well defined morphism of ${\mathscr O}$-modules between the cohomologies, 
$\pmb\xi^j:{\mathscr H}^j({\mathscr F})\to{\mathscr H}^j({\mathscr J})$.
On the other hand, as $\phi$ and $\psi$ are gauge fields, $\hat\pi\hat\xi=0$. Hence, $\pmb\xi^j$ defines morphisms of ${\mathscr O}$-modules $\zeta^j({\xi}):{\mathscr H}^j({\mathscr F}) \to\Omega^1\otimes_{\mathscr O}{\mathscr H}^j({\mathscr F}).$ We denote by $\vartheta^j$ and $\chi^j$ the  connections on ${\mathscr H}^j({\mathscr F})$ determined by $\hat\psi$ and $\hat\phi$, respectively.
Since $\pmb\xi^j=(\pmb\phi^j-\pmb\eta^j)-(\pmb\psi^j-\pmb\eta^j)$, it follows that $\chi^j=\vartheta^j+{\zeta}^j(\xi)$. In short, $\zeta^j(\xi)$ is the ``variation'' on the connection $\vartheta^j$ induced by the ``variation'' $\xi$ of the gauge field $\psi$.
\end{Rems}



\subsection{Yang-Mills fields  on a brane.}\label{SsYng-Millsonabrane}

The result deduced in the following paragraph gives us a suggestion for the definition of the  Yang-Mills functional over the gauge fields on a brane.


\subsubsection{An Euler-Poincar\'e mapping.}\label{Sss: ConnectionsCS}  
Let ${\mathscr A}$ be a coherent sheaf on the Hodge manifold $Y,$  and $\alpha:{\mathscr A}\to \Omega^1({\mathscr A})$ a
 holomorphic connection on ${\mathscr A}$. Denoting by ${\mathcal S}_{\mathscr A}$ the singularity set of ${\mathscr A}$,   on
 $Y\setminus{\mathcal S}_{\mathscr A}$ we define differential form 
$$\Phi({\mathscr A},\alpha):={\rm tr}\big(K_{\alpha}\wedge\star K_{\alpha} \big)\in 
\Gamma(Y\setminus{\mathcal S}_{\mathscr A},\,\Omega^{\rm top}\big) ,$$ 
  $K_{\alpha}$ being the  curvature of $\alpha$ considered as an  ${\rm End}(A)$-valued $2$-form. 

By $\mathfrak{C}$, we denote the category whose objects are pairs $({\mathscr A},\,\alpha)$. A morphism $f:({\mathscr A},\,\alpha)\to
({\mathscr B},\,\beta)$ is a morphism of coherent sheaves compatible with the connections; i.e. such that $(1\otimes f)\circ\alpha=\beta\circ  f$.
\begin{Prop}\label{P:EulerPoincsre}
If $0\to ({\mathscr A},\,\alpha)\overset{f}{\to} ({\mathscr B},\,\beta)\overset{g}{\to} ({\mathscr C},\,\gamma)\to 0,$
is an exact sequence in $\mathfrak{C}$, then on $Y\setminus {\mathcal S}$
$$\Phi({\mathscr B},\beta)=  \Phi({\mathscr A},\alpha)+\Phi({\mathscr C},\gamma),$$ 
 where   ${\mathcal S}$ is the union of the singularity sets of ${\mathscr A},$ ${\mathscr B},$ and ${\mathscr C}$.
\end{Prop}
{\it Proof.} Let $y_0\in Y\setminus {\mathcal S}$. As the exact sequence splits locally on $Y\setminus {\mathcal S}$, there exists an open neighborhood $U$ of $y_0$ such that $g|_U$, in the sequence of locally free modules   
 $0\to {\mathscr A}|_U\overset{f|_U}{\to} {\mathscr B}|_U\overset{g|_U}{\to} {\mathscr C}|_U\to 0,$
has a right inverse $h$.

Let $a$ be a frame for ${\mathscr A}|_U$, then $\alpha(a)={\sf A }\cdot a$, where ${\sf A}$ is a matrix of $1$-forms on $U$. Furthermore, $a$ can be chosen so that ${\sf A}(y_0)=0$. Similarly, let $c$ be a frame for ${\mathscr C}|_U$, then $\gamma(c)={\sf C}\cdot c$ and we choose $c$ so that ${\sf C}(y_0)=0$.  From the splitting, it follows that $\{f(a),\, h(c)\}$ is a frame for ${\mathscr B}|_U$.
By the compatibility of the connections with $f$ and $g$, 
$$\beta(f(a))=(1\otimes f)(\alpha(a))=(1\otimes f)({\sf A} \cdot a)= {\sf A}\cdot f(a).$$
On the other hand, $\beta(h(c))={\sf R} \cdot f(a)+{\sf S}\cdot h(c) $, with ${\sf R}$ and ${\sf S}$ matrices of $1$-forms. But,
$$ {\sf C}\cdot c=\gamma(c)=\gamma(gh(c))=(1\otimes g)(\beta(h(c))=(1\otimes g)\big({\sf R}\cdot f(a)+{\sf S} \cdot h(c)\big).$$
As $g\circ f=0$ and $g\circ h=1$, it follows that ${\sf C}={\sf S}.$ That is, the matrix of the connection $\beta$ in the frame $\{f(a),\, h(c)\}$  is
\begin{equation}\label{matrix Kb}
{\sf M}:=\begin{pmatrix}
{\sf A} & {\sf R} \\
0 & {\sf C}
\end{pmatrix}
\end{equation}

Since ${\sf A}(y_0)=0$ and ${\sf C}(y_0)=0$, the matrix of $K_{\alpha}(y_0)$, of the curvature of $\alpha$ at the point $y_0$, is
 $d{\sf A}$. Analogous  the matrix of $K_{\gamma}(y_0)$ is $d{\sf C}$. The one of $K_{\beta}(y_0)$ is the exterior derivative of 
(\ref{matrix Kb}), since ${\sf M}\wedge {\sf M}=0$ at $y_0$. Then 
\begin{align}
{\rm tr }\big( K_{\beta}(y_0) \wedge\star K_{\beta}(y_0)  \big)&={\rm tr}( d{\sf A}\wedge\star d{\sf A} ) +
{\rm tr}( d{\sf C}\wedge\star d{\sf C} ) \notag \\  
&={\rm tr }\big( K_{\alpha}(y_0) \wedge\star K_{\alpha}(y_0)  \big)+{\rm tr }\big( K_{\gamma}(y_0) \wedge\star K_{\gamma}(y_0)  \big).\notag
\end{align}
As $y_0$ is an arbitrary point of $Y\setminus{\mathcal S}$, it follows the proposition.
\qed

\smallskip

Let   $({\mathscr G}^{\bullet},\,{\delta}^{\bullet})$ be a bounded complex of {coherent}  sheaves on the manifold $Y$.
Let $\nabla^{\bullet}$ be a family of holomorphic connections, compatible with the operators $\delta^{\bullet}$. 
That is, $\nabla^i:{\mathscr G}^{i}\rightarrow \Omega^1({\mathscr G}^{i})$ is a holomorphic gauge field on the coherent sheaf ${\mathscr G}^i$ such that
 $$(1\otimes \delta^i)\nabla^i=\nabla^{i+1} \delta^i.$$
Hence, $\nabla^i({\mathscr Ker}({\delta}^i))\subset {\mathscr Ker}(1\otimes{\delta}^i)$ and a similar relation for the image
  ${\mathscr Im}(\delta^{i-1})$.
It follows that $\nabla^i$ induces a connection $\theta^i$ on the cohomology
$$\theta^i:{\mathscr H}^i({\mathscr G}^{\bullet})\to {\mathscr H}^i(\Omega^1({\mathscr G}^{\bullet})).$$
 Obviously, the restrictions of $\nabla^i$  determine connections on ${\mathscr Ker}(\delta^i)$ and
 ${\mathscr Im}(\delta^{i+1})$, respectively. one has the exact sequence
 \begin{equation}\label{exact_seq:C}
 0\to \big({\mathscr Ker}(\delta^i),\,\nabla^i\big)\to\big({\mathscr G}^i,\,\nabla^i\big)\to \big({\mathscr Im}(\delta^{i}),\,\nabla^{i+1}\big)\to 0
 \end{equation}
 in the category $\mathfrak{C}$. Similarly, we have the exact sequence
 \begin{equation}\label{exact_seq:C1}
 0\to  \big({\mathscr Im}(\delta^{i-1}),\,\nabla^{i}\big)\to\big({\mathscr Ker}(\delta^i),\,\nabla^i\big)\to\big({\mathscr H}^i,\,\theta^i\big) \to 0.
 \end{equation}

\begin{Cor}\label{Coro:Euler}
Denoting with ${\mathcal S}$ the union of the singularity sets of the sheaves ${\mathscr G}^i$, then on $Y\setminus {\mathcal S}$
 \begin{equation}\label{equalityforms}
\sum_{i}(-1)^i{\rm tr}\big(K_{\nabla^i}\wedge\star K_{\nabla^i}\big)=
  \sum_{i}(-1)^i{\rm tr}\big(K_{\theta^i}\wedge\star K_{\theta^i}\big).
 \end{equation}
\end{Cor}
  {\it Proof.} From Proposition \ref{P:EulerPoincsre} together with (\ref{exact_seq:C}), it follows
		$$\Phi({\mathscr G}^i,\,\nabla^i)= \Phi({\mathscr Ker}(\delta^i),\,\nabla^i)+\Phi({\mathscr Im}(\delta^i),\,\nabla^{i+1}).$$
	From (\ref{exact_seq:C1}), one obtains an analogous relation. Taking the alternate sums 
  $$\sum_i(-1)^i\Phi({\mathscr G}^{i},\,\nabla^{i})=\sum_i(-1)^i\Phi({\mathscr H}^{i},\,\theta^{i}). $$
  \qed

  
	\subsubsection{The Yang-Mills functional}\label{Ss:Yang-Millsfunct}  
We propose a definition for the Yang-Mills functional over gauge fields on a brane. This proposal is based on the following considerations:
\begin{enumerate}
  \item It is reasonable to require that this definition generalizes the one for coherent sheaves.
  \item As a gauge field is a homotopy class of a morphism of complexes, it seems convenient to move on  the cohomology of these complexes.
  \item Let $E^{\bullet}$ be a bounded complex of   Hermitian vector bundles over the K\"ahler manifold $Y,$ and $\nabla^{\bullet}$ a family of connections compatible with 
	the Hermitian metrics and  
	the coboundary operators. Denoting by $H^i(E^{\bullet})$ the cohomology bundles, there exist connections $\theta^i$ on those bundles, induced by the family $\nabla^{\bullet}$. By Corollary \ref{Coro:Euler} together with   (\ref{traza-Norma}), one has
  the following equality of Euler-Poincar\'e type. 
  $$\sum_i(-1)^i\norm{K_{\nabla^i}}^2  =
  \sum_i(-1)^i\norm{(K_{\theta^i}}^2.$$ 
 \end{enumerate}
 On the basis of the above considerations, it seems appropriate to define the value of the Yang-Mills functional on the    gauge $\psi$ on the brane ${\mathscr F}^{\bullet}$ as
 $\sum_i(-1)^i\norm{{\mathcal K}_{\vartheta_i}}^2.$ 
 
   More precisely, taking into account Proposition \ref{P:Connvartheta}, we adopt the following definitions.   

\begin{Def}\label{D:HermitianBrane}
The brane $({\mathscr F}^{\bullet},\,\delta^{\bullet})$ is called a Hermitian brane, if the cohomology sheaves ${\mathscr H}^j$ are Hermitian.
\end{Def}

Let $({\mathscr F}^{\bullet},\,{\delta}^{\bullet})$ be a Hermitian brane on the Hodge manifold $Y$. Given a gauge
 field $\psi$  
on the brane $({\mathscr F}^{\bullet},\,\delta^{\bullet})$,
by Proposition \ref{P:Connvartheta}, one has the family of curvatures ${\mathcal K}_{\vartheta^i}$ of the connections induced on the cohomologies.
	We denote by ${\mathcal S}^i$ the singularity set of the cohomoloy sheaf ${\mathscr H}^i$ and let ${\mathcal S}:=\cup{\mathcal S}^{i}.$ On
	$Y\setminus {\mathcal S} $ all the ${\mathscr O}$-modules ${\mathscr H}^i$ are locally free and we denote by $H^i$ the corresponding vector bundles. One has  the respective curvature $2$-forms
	$$  K_{\vartheta^{\bullet}} \in\Gamma( Y\setminus {\mathcal S},\, \Omega^2\otimes_{\mathscr O}{\mathscr End}({ H}^{\bullet})).$$

\begin{Def}
 Given a gauge field $\psi$ on the Hermitian $B$-brane $({\mathscr F}^{\bullet},\,\delta^{\bullet}),$  we define the value of the Yang-Mills functional at $\psi$ by
 \begin{equation}\label{Y-M(Brana)}
 \mathcal{YM}(\psi)=\sum_i(-1)^i\norm{K_{\vartheta_i}}^2.
 \end{equation}
\end{Def}
	
	The {\it Yang-Mills fields} on the brane 
	$({\mathscr F}^{\bullet},\,{\delta}^{\bullet})$ are 
	 the critical points of this functional.

	Note that, if $({\mathscr F}^{\bullet},\,\delta^{\bullet})$ is an acyclic complex, then the Yang-Mills functional for this complex is identically zero.

		
\Exampl	
Let ${\mathscr A}^{\bullet}:=\big( {\mathscr A}^{\bullet},\,d_A^{\bullet},\,\alpha^ {\bullet }\big)$ be a complex in the category ${\mathfrak C}$; i.e, a complex of coherent sheaves with a family of holomorphic connections compatible with the coboundary operator $d_A$. Let $f:=(f^{\bullet})$ a morphism 
  $f^{\bullet}:{\mathscr A}^{\bullet}\to{\mathscr B}^{\bullet}$ in   ${\mathfrak C}$; that is, $f$ is a morphism of complexes compatible with the connections. One can consider the mapping cone ${\mathscr C}^{\bullet}$ of 
  $f$. Thus, 
${\mathscr C}^{\bullet}=\big( {\mathscr A}^{\bullet}[1] \oplus {\mathscr B}^{\bullet},d_C^{\bullet},\nabla^{\bullet},
 \big)$, with $d_C(a,\,b)=\big(d(a),\,(-1)^{{\rm degree}\,a}f(a)+db   \big)$ and $\nabla(a,\,b)=(\alpha(a),\,\beta(b)).$ In fact, $(1\otimes d_C)\circ\nabla=\nabla\circ d_C$ and thus ${\mathscr C}^{\bullet}$ is a complex of the category ${\mathfrak C}$
 
 For each $i$ one has the following exact sequence in the category ${\mathfrak C}$
 $$0\to{\mathscr B}^{i}\to {\mathscr C}^{i}\to {\mathscr A}^{i}[1]\to 0.$$
 From Proposition \ref{P:EulerPoincsre},
 $\Phi({\mathscr B}^{i})+ \Phi({\mathscr A}^{i+1})=\Phi({\mathscr C}^{i})$. Multiplying by $(-1)^i$ and summing
\begin{align} \label{alingsumi}\sum_i(-1)^i{\rm tr}(K_{\beta^i}\wedge\star K_{\beta^i}) + &
 \sum_i(-1)^i{\rm tr}(K_{\alpha^{i+1}}\wedge\star K_{\alpha^{i+1}}) \\ \notag
 &=\sum_i(-1)^i{\rm tr}(K_{\nabla^i}\wedge\star K_{\nabla^i}).
\end{align}

Let us assume that
\begin{itemize}
 \item ${\mathscr A}^{i}$ and ${\mathscr B}^{i}$   Hermitian sheaves for all $i$.
 \item  $\alpha^i$ and $\beta^i$  are Hermitian gauge fields (i.e., compatible with the metric) on ${\mathscr A}^{i}$ and ${\mathscr B}^{i}$, respectively.
 \end{itemize}
Then one defines on ${\mathscr C}^i$  the metric
$\langle(a,b),\,(a',b')\rangle:=\langle a,\,a'\rangle+\langle b,\,b'\rangle$. The connection $\nabla^i$ is compatible with this metric.
 From the  equality  (\ref{alingsumi})  together with (\ref{traza-Norma}), one deduces the following proposition. 
\begin{Prop}\label{P:YMCone}
 With the above notations and under the above hypotheses,  
  $\alpha$ and $\beta$ determine in a natural way a gauge field $\nabla$ on the mapping cone of $f^{\bullet}$ satisfying 
 \begin{equation}\label{YMCone}
  \mathcal{YM}(\beta)-\mathcal{YM}(\alpha) = \mathcal{YM}(\nabla).
  \end{equation}
\end{Prop}
On the other hand, in the context of the branes theory,
the fact that the branes ${\mathscr  A}^{\bullet},\,{\mathscr  B}^{\bullet}$ and ${\mathscr  C}^{\bullet}$  are the members of 
the distinguished triangle ${\mathscr A}^{\bullet}\to {\mathscr B}^{\bullet}\to 
{\mathscr C}^{\bullet}\overset{+1}{\to}$ means that 
${\mathscr A}^{\bullet}$ and  ${\mathscr C}^{\bullet}$   can potentially bind together to form the  brane ${\mathscr B}^{\bullet}$ \cite[Section 6.2.1]{Aspin}. Thus, the additive nature of equation
(\ref{YMCone}) is consistent with this interpretation.


\subsubsection{Yang-Mills fields}\label{Ss:Y_Mbrane}
	If $\phi$ and $\psi$ are
	gauge fields on the $B$-brane $({\mathscr F}^{\bullet},\,\delta^{\bullet})$ and $\xi=\phi-\psi$, using the notations introduced in Remark \ref{R:variation},
	the connections on the cohomologies induced by $\phi$ and $\psi$ satisfy $\chi^j(\xi)=\vartheta^j+\zeta^j(\xi)$, with 
	$$\zeta^j(\xi)\in \Gamma(Y\setminus {\mathcal S},\,\Omega^1\otimes_{\mathscr O}{\mathscr End}(H^j)).$$
	
	With the mentioned notation, 
	an infinitesimal  variation $\psi_{\epsilon}$ of $\psi$ is given by a family 
	$\epsilon^j\hat\xi^j$, with $\epsilon^j\in{\mathbb C}$,
	which defines a morphism between ${\mathscr F}^{\bullet}$ and 
	$\widehat{{\mathscr J}({\mathscr F}^{\bullet})}$ in the homotopy category. In this case, for the  connections on the cohomologies, one has  $\chi^j=\vartheta^j+\epsilon^j\zeta^j.$ 
	Furthermore, on $Y\setminus {\mathcal S}$ the curvatures satisfy
	$$K_{\chi^j}=K_{\vartheta^j}+\epsilon^j\vartheta^j({\zeta^j})+O((\epsilon^j)^2),$$
	$\vartheta^j({\zeta^j})$ being the covariant derivative of ${\zeta^j}$ considered as a section 
	of $\Omega^1\otimes_{\mathscr O}{\mathscr End}(H^j)$.
	
	The functional $\mathcal{YM}$ takes at the  gauge field ${\psi}$ a stationary value if
	$$\frac{\partial}{\partial\epsilon^i}\mathcal{YM}(\psi_{\epsilon} )\Big|_{\epsilon^{\bullet}=0}=0,$$
	for all $i$ and any  variation of $\psi$. That is, if 
	\begin{equation}\label{langleK_}
	\langle K_{\vartheta^i},\,\vartheta^i(\zeta^i)\rangle=0,
	\end{equation}
   for all $i$ and for any $\zeta^i$ defined by a  variation of $\psi$. 
	Therefore, by (\ref{YanMillsEq}),
	one has the following proposition.
	\begin{Prop}\label{estationaryY-Mfunctional}
Let	$({\mathscr F}^{\bullet},\,\delta^{\bullet})$
	be a Hermitian brane,  such that the sheaves ${\mathscr H}^i({\mathscr F}^{\bullet})$ are reflexive, and let $\psi$  gauge field   on ${\mathscr F}^{\bullet}$.
	If  $\vartheta^i$ is a Yang-Mills field
	on ${\mathscr H}^i$ for all $i$, then $\psi$
	is  a stationary point  of the Yang-Mills functional; i.e. $\psi$ is a Yang-Mills field on the brane.
	\end{Prop}

	The following proposition is a  converse  to  Proposition \ref{estationaryY-Mfunctional}. 
	\begin{Prop}\label{Th:inverse} 
	 Let ${\mathscr F}^{\bullet}$ be a Hermitian  $B$-brane as in Proposition \ref{estationaryY-Mfunctional}. 
	If $\nabla^{\bullet}$ is a Yang-Mills field on   ${\mathscr F}^{\bullet}$, then the connection $\vartheta^j$ induced on ${\mathscr H}^j$ is a Yang-Mills field on this sheaf. 
	\end{Prop}
 {\it Proof.} We will consider $({\mathscr F}^{\bullet},\,\delta^{\bullet})$ as an object of the category
 $K^b_{\rm coh}(Y)$; that is, we assume that  $({\mathscr F}^{\bullet},\,\delta^{\bullet})$ is a complex of injective ${\mathscr O}$-modules such that its cohomology is bounded and coherent.
 
 As ${\mathscr Ker}(\delta^i)$ is a submodule of the injective ${\mathscr O}$-module ${\mathscr F}^i$, then ${\mathscr Ker}(\delta^i)$ is a retract of ${\mathscr F}^i$ and thus it is also an injective ${\mathscr O}$-module. Therefore,   the following short exact sequence
 $$0 \to {\mathscr Ker}(\delta^i) \to {\mathscr F}^i\to {\mathscr Coim}(\delta^i)\to 0$$
 splits. That is, ${\mathscr F}^i\simeq {\mathscr Ker}(\delta^i) \oplus {\mathscr Coim}(\delta^i).$
 
 Since ${\mathscr Im}({\delta}^{i-1})$ is a retract of the injective ${\mathscr O}$-module ${\mathscr Ker}(\delta^i)$, the following ahrt exact sequence also splits
 $$0\to {\mathscr Im}({\delta}^{i-1})\to {\mathscr Ker}(\delta^i) \to {\mathscr H}^i\to 0.$$
 Thus, 
 \begin{equation}\label{Fi=Hi+Gi}
 {\mathscr F}^i\simeq {\mathscr H}^i \oplus {\mathscr G}^i,
 \end{equation}
 where ${\mathscr G}^i$ is isomorphic to the direct sum of ${\mathscr Coim}(\delta^i)$ and ${\mathscr Im}(\delta^{i-1})$.
 As ${\mathscr H}^i$ and ${\mathscr G}^i$ are summands in a direct sum decomposition of an injective ${\mathscr O}$-module, they are also injective.
 
 On the other hand, the coboundary operator $\delta^i:{\mathscr F}^i\to {\mathscr F}^{i+1}$ induces via the isomorphisms  (\ref{Fi=Hi+Gi}) to the morphism
 \begin{equation}\label{Fi=Hi+Gi(1)}
  \delta^i:{\mathscr H}^i \oplus {\mathscr G}^i\to {\mathscr H}^{i+1} \oplus {\mathscr G}^{i+1},\;\;\;
  (a,\,b)\mapsto (0,\,\delta^ib).
  \end{equation}
 
 Given $\xi\in {\rm Hom}_{K^b_{\rm coh}(Y)}\big({\mathscr F}^{\bullet},\,\Omega^1\otimes_{\mathscr O}{\mathscr F}^{\bullet}  \big)$, according to Remark   \ref{R:variation},  it determines $\zeta^i\in {\rm Hom}\big({\mathscr H}^{i},\,\Omega^1\otimes_{\mathscr O}{\mathscr H}^{i}  \big).$  As $\nabla^{\bullet}$ is a Yang-Mills field (\ref{langleK_}) is satisfied. 
 
 A general ``variation'' of $\vartheta^j$ is defined by an element
 $\tau\in{\rm Hom}\big({\mathscr H}^j,\,\Omega^1\otimes_{\mathscr O}{\mathscr H}^j\big)$. We need to prove that
 $$\langle K_{\vartheta^i},\,\vartheta^i(\tau)\rangle=0,$$
 for any variation $\tau$.
 
 The morphism $\tau$ can be extended 
 $$C^i:{\mathscr H}^i\oplus{\mathscr G}^i\to \Omega^1\otimes_{\mathscr O}\big( {\mathscr H}^i\oplus{\mathscr G}^i  \big),$$
 where 
 $$C^i(a,\,b)=\begin{cases}(\tau(a),\,0),\;\;\text{if}\;\, i=j\\
 (0,\,0), \;\;\text{if}\;\, i\ne j
 \end{cases} $$
 Moreover, the $C^i$ are compatible with the coboundaries. For example for $i=j$, by  (\ref{Fi=Hi+Gi(1)}),
 $((1\otimes \delta^j)\circ { C}^j)(a,\,b)=(1\otimes\delta^j)(\tau(a),\,0)=0$; and 	${ C}^{j+1}\circ\delta^{j}(a,\,b)=0$. Thus, by the isomorphism (\ref{Fi=Hi+Gi}) the $C^i$ determine a morphism
 $\xi:{\mathscr F}^{\bullet}\to\Omega^1\otimes_{\mathscr O}{\mathscr F}^{\bullet}$ in the category $K^b_{\rm coh}(Y)$, and the corresponding $\zeta^i$ induced in the cohomologies are all $0$ except when $i=j$, in which case $\zeta^j=\tau$. Hence,
  by (\ref{langleK_}),
	$$\langle K_{{\vartheta}^j},\,\vartheta^j(\tau)\rangle= \langle K_{{\vartheta}^j},\,\vartheta^j(\zeta^j)\rangle=0.$$
	This holds for any ``variation'' $\tau$ of $\vartheta^j$.
	That is, by (\ref{YanMillsEq}), $\vartheta^j$ is a Yang-Mills field on ${\mathscr H}^j$. \qed
\smallskip

\noindent
{\bf Proof of Theorem \ref{C:definitivo}.}	
 From Proposition \ref{Th:inverse} together with Proposition \ref{estationaryY-Mfunctional}, it follows Theorem \ref{C:definitivo}
 \qed
 
	\smallskip
	
Let us assume that the set of gauge fields on the brane ${\mathscr F}^{\bullet}$ is nonempty. Let $m$ be the dimension of the vector space ${\rm Ext}^0({\mathscr F}^{\bullet},\,\Omega^1({\mathscr F}^{\bullet})).$ We denote by $\xi_1,\dots, \xi_m$ a basis of this vector space. According to Proposition 
\ref{P:Atiyah}, any gauge field $\psi$ on the brane can be expressed 
$$\psi=\tilde\psi+\sum_a\lambda_a\xi_a$$
 $\tilde\psi$ being a fixed gauge field and $\lambda_a\in{\mathbb C}.$ The corresponding connection $\nabla^i$ on ${\mathscr F}^i$ is of the form
  $$\nabla^i=\tilde\nabla^i+\sum_a\lambda_a\xi_a^ i,$$
 with $\xi^i_a\in {\rm Hom}({\mathscr F}^i,\,\Omega^1({\mathscr F}^i))$. Hence, the connections on the cohomology sheaves ${\mathscr H}^i$ 
 can be written in the form
 $$\vartheta^i=\tilde\vartheta^i+\sum_ a\lambda_a \zeta^i_a.$$
 The corresponding curvatures satisfy
 $$K_{\vartheta^i}= K_{\tilde\vartheta^i}+\sum_a \lambda_a\tilde\vartheta^i(\zeta^i_a)+\sum_{a,b}\lambda_a\lambda_b\zeta^i_a\wedge\zeta^i_b.$$
 Therefore 
$\norm{K_{\vartheta^i}}^2$ is a   polynomial $P^i(\lambda_1,\dots,\lambda_m)$ of degree $\leq 4.$

By Theorem \ref{C:definitivo}, the critical points of the Yang-Mills functional on ${\mathscr F}^{\bullet}$  correspond to the points $(\lambda_1,\dots,\lambda_m)\in{\mathbb C}^m$ which satisfy the equations
$\frac{\partial P^i}{\partial\lambda _a}=0,$  for $a=1,\dots, m$ and for all $i$. 
 We have the following result, which generalizes Theorem \ref{P:numeroYM}.
\begin{Thm}\label{P:numero YM} 
Assumed the cohomology sheaves ${\mathscr H}^i$ of the brane ${\mathscr F}^{\bullet}$ are reflexive and  the set of gauge fields on   ${\mathscr F}^{\bullet}$ is nonempty. Then this set
is in bijective correspondence with the points of a subvariety of ${\mathbb C}^m$ defined by $m\cdot s$  polynomials of degree $\leq 3$,  where  
$m:={\rm dim}\,{\rm Ext}^0({\mathscr F}^{\bullet},\,\Omega^1({\mathscr F}^{\bullet}))$ and   $s$ the number of nontrivial
sheaves ${\mathscr H}^i$. 
	\end{Thm} 
	
	
	\section{Appendix}
	Let ${ L}$ be a holomorphic line bundle over the K\"ahler manifold $Y$. We denote by ${\mathcal U}=\{U_a\}_a$
	a good cover of $Y,$ such that the restrictions ${ L}|_{U_a}$ are trivial. Let $\{\varphi_{ab}\}$ denote the corresponding cocycle of ${\mathcal Z}^1({\mathcal U},\,{\mathscr O}^{\times})$. By the simply connectedness of each $U_{ab}$, one can define  $\frac{1}{2\pi i}\log \varphi_{ab}$, and
	 $\zeta:=\{\frac{1}{2\pi i}\partial \log \varphi_{ab}\}$ is a $\check{\rm C}$ech cocycle in ${\mathcal Z}^1({\mathcal U},\,{\Omega}^1)$. 
	 \begin{Lem} \label{L:Appendix} The cocycle $\zeta$
	   is a coboundary iff the Atiyah class $a({ L})=0$.
	 \end{Lem}
	{\it Proof.} A holomorphic gauge field  in  the trivialization $L|_{U_a}$  is of the form $\partial +B_a$, with $B_a$ a holomorphic $1$-form on $U_a$.  The local connections $\{\partial+B_a\}$ can be glued to form a holomorphic connection
	on $L$ iff on $U_{ab}$
	$$B_b-B_a=\partial\log \varphi_{ab},$$
	for all $a,b$. Equivalently, if the cocycle $\zeta$ is a coboundary. 
	\qed
	
\smallskip

If $\beta\in{\mathcal Z}^0({\mathcal U},\,\Omega^1)$ satisfies $\delta\beta=\zeta$, where $\delta$ is the $\check{\rm C}$ech  coboundary operator, then $-d\beta$ is the $2$-form on $Y$ determined by $\zeta$. Such a $\beta$ can be construct from a Hermitian metric on $L$.	
	
	A Hermitian metric on $L$ is defined by a family $\{f_a:U_a\to{\mathbb R}_{>0}\}_a$ of ${C}^{\infty}$
	functions, such that $f_b=\varphi_{ab}\bar \varphi_{ab}f_a$. Letting $\beta_a=\frac{1}{2\pi i}\partial\log f_a$,
	since the transition functions $\varphi_{ab}$ are holomorphic one has
	 $$(\delta\beta)_{ab}=\frac{1}{2\pi i}\partial\log(\varphi_{ab}\bar \varphi_{ab})=\zeta_{ab}.$$
	On the other hand,
	 $$-d\beta=\frac{i}{2\pi}\bar\partial\partial\log f_a.$$
	 Thus,
	 $$(\delta(-d\beta))_{ab}=\frac{-1}{2\pi i}\bar\partial\partial\log(\varphi_{ab}\bar \varphi_{ab})=0.$$
	  That is, $-d\beta$ is a  $2$-form on $Y$.
	
	
	   
	 Associated to the Hermitian metric  is defined the corresponding Chern connection, whose curvature form is given by
	 $\bar\partial\partial\log f_a$. Therefore, the cocycle $\zeta$ defines the first $c_1(L)$. From  Lemma \ref{L:Appendix}, it follows the following known result \cite{Atiyah}.
	 
	 \begin{Prop}\label{a=c-1} $L$ admits a holomorphic gauge field iff $c_1(L)=0.$
	 \end{Prop}
	  Thus, $L$ supports a holomorphic gauge field iff $L$ is flat.



\end{document}